# The sequence of Collatz functions, exceptionality of the 3n+1 function and the notion of Collatz generalized Matrix

Raouf Rajab
raouf.rajab@enig.rnu.tn

## Abstract

In the present paper, we define a sequence of functions which has very important proprieties. Each function which belongs of this sequence of functions behaves in a very similar way to the Collatz function and it can generate an integer sequences having the same properties of Collatz sequences. For any integer n, we define a subset of $\mathbb{Z}$ denoted by **D(n)** which contains all **integers greater than or equal to (2n+2)** next, we define the function denoted by $\mathrm{F_n}$ **from D(n) to $\mathbb{Z}$** such that the value of $\mathrm{F_n(P)}$ is determined as follows: if **P** is **even** , we multiply it by 3 then we subtract **2n** from **3P** and the value obtained **(3P-2n)** is divided by 2 so $\mathrm{F_n(P)}$ equal to **one half of (3P-2n),** if **P** is **odd** then $\mathrm{F_n(P)}$ equal to **one half of (P+2n+1)**. The first property obtained is that the pair of integers **(2n+2, 2n+3)** is stable by iteration of the $\mathrm{F_n}$ function, in other words:

$$\mathbf{F_n(2n+2) = 2n+3\ ;\ F_n(2n+3) = 2n+2}$$

An important theorem is also proved in this article which concerned the relation between Collatz function denoted by T and $\mathrm{F_n}$ function. This relation can be expressed by the following equation:

$$\mathbf{T^k(P) = F_n^k(P+2n+1) - (2n+1)} \text{ where } \mathrm{k \in \mathbb{N}; P \in \mathbb{N}^*\ ; n \in \mathbb{Z}}$$

Based on this theorem, we prove that Collatz conjecture can be stated as follows: Start with any integer P such that P is **greater than or equal to (2n+2)** , the conjecture is that no matter what value of P , the sequence which generated by $\mathrm{F_n}$ function will always reach to( 2n+2).

**Résumé :** Dans cet article, on définit une suite des fonctions possédant des propriétés très importantes en effet chaque fonction appartenant à cette suite possède un comportement très voisin du comportement de la fonction de Collatz et peut générer des suites possédant des propriétés similaires aux suites générées par le processus itératif de la fonction de Collatz. Pour un entier relatif n quelconque, on définit une partie de $\mathbb{Z}$ qu' on la note D(n) qui contient tous les entiers relatifs supérieurs ou égaux à 2n+2 et pour tout entier relatif P appartenant à D(n) , on définit la fonction $\mathrm{F_n}$ tel que la valeur de $\mathrm{F_n(P)}$ est déterminée de la manière suivante : Si P est pair on le multiplie par 3 ensuite on fait soustraire 2n de 3P puis la valeur obtenue (3P-2n) est devisée par 2 , si P est impair on l'ajoute 2n+1 et la somme obtenue (P+2n+1) est devisée par 2. La première propriété qu'on peut la vérifier aisément est que le couple (2n+2, 2n+3) est un cycle stable par itération pour la fonction $\mathrm{F_n}$ c'est à dire que :

$$\mathbf{F_n(2n+2) = 2n+3\ ;\ F_n(2n+3) = 2n+2.}$$

Un important théorème est également démontré dans ce travail portant sur la relation entre la fonction de Collatz notée T et la fonction $\mathrm{F_n}$ cette relation est traduite sous forme de l équation suivante:

$$\mathbf{T^k(P) = F_n^k(P+2n+1) - (2n+1);\ Avec\ k \in \mathbb{N}\ ; P \in \mathbb{N}^*\ et\ n \in \mathbb{Z}}$$

En se basant ce théorème, on montre que la conjecture de Collatz peut être énoncée pour chaque fonction $\mathrm{F_n}$ de la manière suivante : Pour tout entier relatif P tel que P est supérieur ou égal à 2n+2, la suite générée par la fonction $\mathrm{F_n}$ finit toujours par atteindre 2n+2 après un nombre fini des itérations.

## 1. Introduction et préliminaire

La fonction de Collatz fait partie d'un ensemble des fonctions numérique particulières. La propriété principale de ce genre des fonctions est que l'expression de l'image d'un entier naturel (ou relatif) dépend de sa parité (ou de sa forme). Cette dépendance entre la forme d'un entier et l'expression de son image engendre un certain nombre des propriétés très exceptionnelles pour les suites de Collatz et aussi pour n'importe quelles suites générées par des fonctions basées sur ce même principe. Plusieurs travaux ont été intéressés par l'étude des fonctions similaires à la fonction des Collatz et par l'étude des suites générées par ces fonctions particulières dans l'objectif est de mieux comprendre les comportements des suites de Collatz et d'expliquer leurs convergence vers le cycle (1,2) ou vers un autre cycle stable par itérations [1],[2],[3] ],[4]. Dans ce travail, on commence par la définition de la fonction de Collatz puis on définit une suite des fonctions qu'on la note $(F_n)_{n\in\mathbb{Z}}$ tel que l'expression de chaque fonction est donnée en fonction d'un entier relatif n quelconque. Des résultats extrêmement importants obtenus sont présentés brièvement dans cette section puis sont étudiés par détailles et démontrés dans le reste de cet article.

**Notation 1.1**

Pour tout entier naturel non nul N, on définit la fonction du Collatz notée T comme suit [1]:

(1.1) $$T(N) = \begin{cases} \frac{N}{2} & \text{si } N \equiv 0 \ (\text{mod}2) \\ \frac{3}{2}N + \frac{1}{2} & \text{si } N \equiv 1 (\text{mod}2) \end{cases}$$

Enoncé de la conjecture de Collatz [1],[2],[4] :

La conjecture de Collatz énonce que $\forall\, N \in \mathbb{N}^*$, il existe un entier naturel non nul k bien déterminé tel que :

$$T^k(N) = 1$$

**Définition 1.1**

Pour un entier relatif bien déterminé n, on définit le sous-ensemble de $\mathbb{Z}$ noté D(n) comme la partie de $\mathbb{Z}$ qui renferme tous les entiers relatifs supérieurs ou égaux à (2n+2). Il est comme suit :

(1.2) $$D(n) = \{P \in \mathbb{Z} \mid P \geq 2n + 2\} = \{2n + 2, 2n + 3, 2n + 4, \dots, 2n + P, \dots\}$$

**Définition 1.2**

Soit n un entier relatif quelconque ($n \in \mathbb{Z}$) , on définit la fonction $\mathbf{F_n}$ sur D(n) tel que pour tout **entier relatif** $P \geq 2n+2$, $F_n(P)$ est exprimée en fonction de n comme suit:

(1.3)
$$\mathbf{F_n(P)} = \begin{cases} \mathbf{\frac{3}{2}P - n} & \textbf{si } \mathbf{P \equiv 0\ (mod2)} \\ \mathbf{\frac{P}{2} + \frac{2n+1}{2}} & \textbf{si } \mathbf{P \equiv 1 (mod2)} \end{cases}$$

**Exemples 1.1**

On donne quelques exemples des fonctions $F_n$ pour les cas ou n=-2, n=-1, n=0,n=1 et n=2 :

$$F_0(P) = \begin{cases} \frac{3}{2}P & \text{si } P \equiv 0 (\text{mod}2) \\ \frac{P}{2} + \frac{1}{2} & \text{si } P \equiv 1\ (\text{mod}2) \end{cases} ;\quad F_1(P) = \begin{cases} \frac{3}{2}P - 1 & \text{si } P \equiv 0\ (\text{mod}2) \\ \frac{P}{2} + \frac{3}{2} & \text{si } P \equiv 1\ (\text{mod}2) \end{cases} ;\quad F_2(N) = \begin{cases} \frac{3}{2}P - 2 & \text{si } P \equiv 0\ (\text{mod}2) \\ \frac{P}{2} + \frac{5}{2} & \text{si } P \equiv 1\ (\text{mod}2) \end{cases}$$

$$F_{-1}(P) = \begin{cases} \frac{3}{2}P + 1 & \text{si } P \equiv 0\ (\text{mod}2) \\ \frac{P}{2} - \frac{1}{2} & \text{si } P \equiv 1 (\text{mod}2) \end{cases} ;\quad F_{-2}(P) = \begin{cases} \frac{3}{2}P + 2 & \text{si } P \equiv 0\ (\text{mod}2) \\ \frac{P}{2} - \frac{3}{2} & \text{si } P \equiv 1 (\text{mod}2) \end{cases}$$

**Les principaux résultats prouvés dans cet article**

(1) Le cycle (2n+2,2n+3) est un cycle stable par itérations pour la fonction $F_n$ avec n un entier relatif quelconque ($n \in \mathbb{Z}$) c'est à dire que :

$$\begin{cases} \mathbf{F_n(2n+2) = 2n+3} \\ \mathbf{F_n(2n+3) = 2n+2} \end{cases}$$

(2) La relation entre T et $F_n$ est décrite par l'équation suivante :

$$\mathbf{T^k(N) = F_n^k(N + 2n + 1) - (2n + 1)}$$

Avec :

$$k \in \mathbb{N}; N \in \mathbb{N}^* \text{ et } n \in \mathbb{Z}$$

(3 )En se basant sur cette dernière équation, on peut montrer que la conjecture de Collatz peut être énoncée pour chaque fonction $F_n$ de la manière suivante :

**$\forall\, P \in \mathbb{Z}$ tel que $P \geq 2n+2$, il existe un entier naturel non nul k tel que :**

$$\mathbf{F_n^k(P) = 2n + 2}$$

Autrement les suites générées par la fonction $F_n$ et dont les premiers termes appartenant à D(n) se comportement de la même manière que les suites de Collatz dont les premiers termes sont des entiers naturels non nuls.

(4) Un autre important théorème est démontré dans ce travail et qui porte sur la nature de la fonction de Collatz T. Ce théorème énonce que pour tous entiers naturels non nuls k et N on peur écrire :

$$T^k(N) = \lim_{n\to+\infty}\left(\frac{1}{2(n+1)}\left(\sum_{j=0}^{n}\left(\left(F_j^k\left(N+A_j\right)+F_{-j-1}^k\left(N+A_{-j-1}\right)\right)\right)\right)\right)$$

Avec: $\mathbf{A_j = 2j+1}$ ; $\mathbf{A_{-j} = -2j-1}$

Ce qui signifie que chaque terme $T^k(N)$ de la suite de Collatz peut être interprété comme la moyenne arithmétique de tous les termes $(F_n^k(N+(2n+1)))_{n\in\mathbb{Z}}$

(5) On sait que n ne peut prendre que des entiers relatifs mais si on fait une seule exception et on prend $\mathbf{n = -\frac{1}{2}}$ on peut tirer deux conséquences importantes:

Si on remplace n par $-\frac{1}{2}$ dans la matrice représentée sur la figure 3 (matrice contenant les suites générées par la fonction $F_n$) on obtient la matrice de Collatz. on peut prendre un exemple

$$F_n(2n+P_i) = 2n+P_{i+1} \text{ pour } \forall\ i \geq 1$$

| | | | | | | | | | |
|---|---|---|---|---|---|---|---|---|---|
| $2n+P_1$ | $2n+P_2$ | $2n+P_3$ | | | $2n+P_i$ | $2n+P_{i+1}$ | | | $2n+P_k$ |

↓ $n = -\frac{1}{2}$

| | | | | | | | | | |
|---|---|---|---|---|---|---|---|---|---|
| $P_1-1$ | $P_2-1$ | $P_3-1$ | | | $P_i-1$ | $P_{i+1}-1$ | | | $P_k-1$ |

$$T(P_i-1) = P_{i+1}-1 \text{ pour } \forall\ i \geq 1$$

Figure1 : Déduction de termes de la suite de Collatz à partir de la suite créée par une fonction $F_n$ quelconque

**Exemple 1.2**

Prenons l'exemple suivant de la suite de premier terme 2n+16 générée par la fonction $F_n$ :

| | | | | | | | | | |
|---|---|---|---|---|---|---|---|---|---|
| $F_n^k(2n+16)$ | 2n+24 | 2n+36 | 2n+54 | 2n+81 | 2n+41 | 2n+21 | 2n+11 | 2n+6 | 2n+9 |

↓ $n = -\frac{1}{2}$

| | | | | | | | | | |
|---|---|---|---|---|---|---|---|---|---|
| 16-1=15 | 24-1=23 | 36-1=35 | 54-1=53 | 81-1=80 | 40 | 20 | 10 | 5 | 8 |

La première ligne contient la suite $S_{F,n}(2n+16,k)$ les termes sont exprimés en fonction de n pour n un entier relatif quelconque. lorsque on remplace n par $(-0{,}5)$, on obtient la suite de Collatz de premier terme 15.

Concernant la matrice généralisée de Collatz, sur la figure suivante, on représente deux matrices dont la première contient des suites générées par la fonction $F_n$ et la deuxième est une matrice qui contient des suites de Collatz. Les cases sont colorées selon les règles suivantes :

| | | |
|---|---|---|
| Suites de Collatz | 2j | 2j + 1 |
| Suites générées par $F_n$ | 2n + (2j) | 2n + (2j + 1) |

Figure 2: Les règles de la coloration de différentes cases pour les différentes suites considérées.

Alors dans le tableau qui contient les suites de Collatz, une case qui contient un entier **impair** sera colorée en bleu alors pour un tableau qui contient les suites générées par la fonction $F_n$ , les cases qui contiennent des entiers **pairs** sont celles qui sont colorées en bleu.

| $P_i$ | $F_n^1(P_i)$ | $F_n^2(P_i)$ | $F_n^3(P_i)$ | $F_n^4(P_i)$ | $F_n^5(P_i)$ |
|---|---|---|---|---|---|
| | | | | | |
| | | | | | |
| | | | | | |
| 2n + 16 | 2n + 24 | 2n + 36 | 2n + 54 | 2n + 81 | 2n + 41 |
| 2n + 15 | 2n + 8 | 2n + 12 | 2n + 18 | 2n + 27 | 2n + 14 |
| 2n + 14 | 2n + 21 | 2n + 11 | 2n + 6 | 2n + 9 | 2n + 5 |
| 2n + 13 | 2n + 7 | 2n + 4 | 2n + 6 | 2n + 9 | 2n + 5 |
| 2n + 12 | 2n + 18 | 2n + 27 | 2n + 14 | 2n + 21 | 2n + 11 |
| 2n + 11 | 2n + 6 | 2n + 9 | 2n + 5 | 2n + 3 | 2n + 2 |
| 2n + 10 | 2n + 15 | 2n + 8 | 2n + 12 | 2n + 18 | 2n + 27 |
| 2n + 9 | 2n + 5 | 2n + 3 | 2n + 12 | 2n + 3 | 2n + 2 |
| 2n + 8 | 2n + 12 | 2n + 18 | 2n + 27 | 2n + 14 | 2n + 21 |
| 2n + 7 | 2n + 4 | 2n + 6 | 2n + 9 | 2n + 5 | 2n + 3 |
| 2n + 6 | 2n + 9 | 2n + 5 | 2n + 3 | 2n + 2 | 2n + 3 |
| 2n + 5 | 2n + 3 | 2n + 2 | 2n + 3 | 2n + 2 | 2n + 3 |
| 2n + 4 | 2n + 6 | 2n+9 | 2n + 5 | 2n + 3 | 2n + 2 |
| 2n + 3 | 2n + 2 | 2n + 3 | 2n + 2 | 2n + 3 | 2n + 2 |
| 2n + 2 | 2n + 3 | 2n + 2 | 2n + 3 | 2n + 2 | 2n + 3 |

| $N_i$ | $T^1(N_i)$ | $T^2(N_i)$ | $T^3(N_i)$ | $T^4(N_i)$ | $T^5(N_i)$ |
|---|---|---|---|---|---|
| | | | | | |
| | | | | | |
| | | | | | |
| 15 | 23 | 35 | 53 | 80 | 40 |
| 14 | 7 | 11 | 17 | 26 | 13 |
| 13 | 20 | 10 | 5 | 8 | 4 |
| 12 | 6 | 3 | 5 | 8 | 4 |
| 11 | 17 | 26 | 13 | 20 | 10 |
| 10 | 5 | 8 | 4 | 2 | 1 |
| 9 | 14 | 7 | 11 | 17 | 26 |
| 8 | 4 | 2 | 1 | 2 | 1 |
| 7 | 11 | 17 | 26 | 13 | 20 |
| 6 | 3 | 5 | 8 | 4 | 2 |
| 5 | 8 | 4 | 2 | 1 | 2 |
| 4 | 2 | 1 | 2 | 1 | 2 |
| 3 | 5 | 8 | 4 | 2 | 1 |
| 2 | 1 | 2 | 1 | 2 | 1 |
| 1 | 2 | 1 | 2 | 1 | 1 |

Figure 3: Similarité entre les comportements des suites $(F_n^k(P_i))_{k\geq 0, P_i\geq 2n+2}$ et des suites de Collatz $(T^k(N_i))_{k\geq 0, N_i\geq 1}$

Ces deux tableaux sont étudiés avec plus des détailles et sont interprétés dans les parties consacrées pour l'étude de différentes fonctions présentées dans cet article.

Les suites de Collatz $Sy(N-1,k)$ peuvent être déduites à partir des suites générées par $F_n$ notées $S_{F,n}(N+2n)$ on fait remplacer n par $-\frac{1}{2}$ dans les expressions de différentes termes de cette dernière suite.

**Les termes impairs dans la matrice de Collatz possède la même distribution des termes pairs dans la matrice globale qui contient les suites générées par $\mathbf{F_n}$.**

**Notations 1.2**

On adopte les notations suivantes :

-Pour une suite de Collatz de premier terme un entier naturel non nul $P_0$ et de longueur finie k :

(1.4) $$Sy(P_0,k) = (P_0, T(P_0), T^2(P_0), T^3(P_0), \dots, T^{k-1}(P_0))$$

-Pour une suite finie générée par la fonction $F_n$ de premier terme un entier relatif $P_0 \geq 2n+2$ et de longueur k, on la note comme suit :

(1.5) $$S_{F,n}(P_0,k) = (P_0, F_n^1(P_0), F_n^2(P_0), F_n^3(P_0), \dots, F_n^{k-1}(P_0))$$

Si les suites considérées ont des longueurs infinies, on les notes comme suit :

$$Sy^{\infty}(P_0) = (P_0, T(P_0), T^2(P_0), T^3(P_0), \dots, T^k(P_0), \dots)$$

$$S_{F,n}^{\infty}(P_0) = (P_0, F_n^1(P_0), F_n^2(P_0), F_n^3(P_0), \dots, F_n^k(P_0), \dots)$$

**Définition 1.3**

Soit N un entier naturel non nul donc pour tout entier naturel n, on définit l'indicateur de parité d'ordre k de P qu'on le note $\mathbf{i}_k(P)$ comme suit :

(1.6) $$\mathbf{i}_k(N) = \begin{cases} 1 \text{ si } T^k(N) \text{ est } \textbf{impair} \\ 0 \text{ si non} \end{cases}$$

De même, on définit l'indicateur de parité relatif à la suite générée par la fonction $F_n$ comme suit :

(1.7) $$\mathbf{f}_{n,k}(P) = \begin{cases} 1 \text{ si } F_n^k(P) \text{ est } \textbf{pair} \\ 0 \text{ si non} \end{cases}$$

Avec $P \in D(n)$

**Notations 1.3**

Le nombre des entiers **impairs** dans une suite de Collatz $Sy(N,k)$ est donné par la relation suivante :

(1.8) $$\beta_k(N) = \sum_{j=0}^{k-1} i_j(N)$$

Le nombre des entiers **pairs** dans une suite $S_{F,n}(P,k)$ est donné par la relation suivante :

(1.9) $$\alpha_{n,k}(P) = \sum_{j=0}^{k-1} f_{n,j}(P)$$

**Notations 1.4**

Soit N un entier naturel non nul et soit P un entier relatif appartenant à D(n) alors on adopte les notations suivantes pour les différentes expressions de $T^k(N)$ et de $F_n^k(P)$ :

-L'expression de $T^k(N)$ en fonction de N et de k est comme suit :

(1.10) $$T^k(N) = \frac{3^{\beta_k(N)}}{2^k} N + r_k(N)$$

Avec :

$\beta_k(N)$ est le nombre des entiers **impairs** dans la suite $Sy(N,k)$

-L'expression de $F_n^k(P)$ en fonction de P et de k est comme suit :

(1.11) $$F_n^k(P) = \frac{3^{\alpha_{n,k}(P)}}{2^k} P + \varphi_{n,k}(P)$$

Avec :

$\alpha_{n,k}(P)$ est le nombre des entiers **pairs** dans la suite $S_{F,n}(P,k)$

Les coefficients $r_k(N)$, $\varphi_{n,k}(P)$ sont appelés les coefficients d'ajustement entière puisqu' ils font ajuster respectivement les valeurs de $T^k(N)$ et de $F_n^k(P)$ à des entiers relatifs.

On ne s'intéresse pas aux expressions de ces coefficients mais on s'intéresse à la détermination de la relation entre $r_k(N)$ et $\varphi_{n,k}(P)$.

**Notation 1.5**

Le vecteur de parité de la suite $Sy(N,k)$ est noté comme suit:

(1.12) $$u(N,k) = \left(\mathbf{i}_0(N), \mathbf{i}_1(N), \mathbf{i}_2(N), \dots, \mathbf{i}_{k-1}(N)\right)$$

Pour une suite de longueur k générée par l'une des fonctions $F_n$, le vecteur de parité est noté comme suit :

(1.13) $$v_n(P,k) = \left(f_{n,0}(P), f_{n,1}(P), f_{n,2}(P), \dots, f_{n,k-1}(P)\right)$$

## 2. Détermination des propriétés des fonctions $(\mathbf{F_n})_{n\in\mathbb{Z}}$

**Définition 2.1**

Pour chaque fonction $F_n$ avec n un entier relatif quelconque, on définit trois points particuliers ou caractéristiques qu'on les notes comme ci-dessous :

(2.1) $$\begin{cases} P_1(n) = 2n + 2 \\ P_2(n) = 2n + 3 \\ A_n = 2n + 1 \end{cases}$$

On peut adopter les notations suivantes pour ces mêmes points caractéristiques.

(2.2) $$\begin{cases} P_1(n) = P_{n,1} \\ P_2(n) = P_{n,2} \end{cases}$$

**Proposition 2.1**

Pour tout **entier relatif** n, le couple $(P_{n,1}, P_{n,2})$ est un cycle stable par itération pour la fonction $F_n$ autrement dit $P_{n,1}$et $P_{n,2}$ satisfaisant les équations suivantes :

(2.3) $$\begin{cases} F_n(2n + 2) = 2n + 3 \\ F_n(2n + 3) = 2n + 2 \end{cases}$$

**Démonstration**

On vérifie que $F_n\big(P_1(n)\big) = P_2(n)$. On écrit alors:

$$\begin{aligned} F_n\big(P_1(n)\big) &= F_n(2n + 2) \\ &= \frac{3}{2}(2n + 2) - n \\ &= 3n + 3 - n \\ &= 2n + 3 \\ &= P_2(n) \end{aligned}$$

De même on peut vérifier que $F_n\big(P_2(n)\big) = P_1(n)$

$$\begin{aligned} F_n\big(P_2(n)\big) &= F_n(2n + 3) \\ &= \frac{2n + 3}{2} + \frac{2n + 1}{2} \\ &= 2n + 2 \\ &= P_1(n) \end{aligned}$$

On peut conclure que :

$$\begin{cases} F_n\big(P_{n,1}\big) = P_{n,2} \\ F_n\big(P_{n,2}\big) = P_{n,1} \end{cases}$$

**Notations 2.2**

Pour un entier relatif n quelconque, le couple constitué par les deux entiers 2n+2 et 2n+3 représente un cycle stable par itération pour la fonction $F_n$ comme on a montré. Il est noté comme ci-dessous:

(2.4) $$C_{F,n} = (2n + 2{,}2n + 3)$$

L'ensemble constitué par les couples qui s'écrivent sous la forme $(2n + 2, 2n + 3)$ pour tout entier relatif n est noté comme :

(2.5) $$\begin{aligned} Cyc(F) &= \{(2n + 2{,}2n + 3) \mid n \in \mathbb{Z}\} \\ &= \{\dots, (-6, -5), (-4, -3)(-2, -1), (0{,}1), (2{,}3), (4{,}5), (6{,}7), (8{,}9), \dots\} \end{aligned}$$

**Remarque 2.1**

Si on remplace n par $-\frac{1}{2}$ dans l'expression des cycles (2n+2,2n+3) on obtient le cycle (1,2) qui n'appartient pas à l'ensemble Cyc(F) et qui représente le cycle stable par itérations pour la fonction de Collatz T de plus ce cycle située entre les deux cycles (0,1) et (2,3) qui correspond aux cas n=-1 et n=0 respectivement.

$-\frac{1}{2}$

| n | -4 | -3 | -2 | -1 | 0 | 1 | 2 | 3 |
|---|---|---|---|---|---|---|---|---|
| $(2n + 2{,}2n + 3)$ | $(-6, -5)$ | $(-4, -3)$ | $(-2, -1)$ | (0,1) | (2,3) | (4,5) | (6,7) | (8,9) |

(1,2)

Figure 4: Exemple illustrant l'exceptionnalité de cycle stable par itérations de la fonction de Collatz

La fonction T de Collatz se comporte comme une fonction très exceptionnelle appartenant à une suite des fonctions telles que les expressions de chaque fonction sont données en fonction d'un entier relatif n bien déterminé et la seule fonction qui peut générer des suites dont les termes sont déduits à partir des suites générées par les fonctions $F_n$en remplaçant n par un entier rationnel négatif.

**Théorème 2.1**

Soient n un entier relatif et k un entier naturel non nul alors pour tous entiers naturels non nuls N, on a:

(2.6) $$T^k(N) + A_n = F_n^k(N + A_n)$$

Cette équation est équivalente à l'équation suivante pour tous entiers relatifs $P \in D(n)$ :

(2.7) $$T^k(P - A_n) = F_n^k(P) - A_n$$

Avec :

$$A_n = 2n + 1$$

**Démonstration**

On montre par récurrence que cette propriété est vraie pour tout entier naturel non nul k. on vérifie tout abord que la propriété est vraie pour i=1.

-Si $N = 2j$

$$\begin{aligned}
F_n^1(N + A_n) - A_n &= F_n^1(2j + 2n + 1) - (2n + 1)\\
&= \frac{2j + 2n + 1}{2} + \frac{(2n + 1)}{2} - (2n + 1)\\
&= \frac{2j}{2} + \frac{(2n + 1)}{2} + \frac{(2n + 1)}{2} - (2n + 1)\\
&= \frac{2j}{2}\\
&= T(2j)\\
&= T(N)
\end{aligned}$$

-Si $N = 2j + 1$

$$\begin{aligned}
F_n^1(N + A_n) - A_n &= F_n^1(2j + 1 + 2n + 1) - (2n + 1)\\
&= F_j^1(2j + 2n + 2) - (2n + 1)\\
&= \frac{3(2j + 2n + 2)}{2} - n - (2n + 1)\\
&= 3j + 3n + 3 - n - 2n - 1\\
&= 3j + 2\\
&= \frac{3}{2}(2j + 1) + \frac{1}{2}\\
&= T(2j + 1)\\
&= T(N)
\end{aligned}$$

Dans les deux cas, on obtient la même expression traduisant la relation entre T et $F_n$ qui s'écrit comme suit:

$$T^1(N) = F_n^1(N + A_n) - A_n$$

La propriété est vraie pour i=1 et on suppose que pour tout entier naturel i allant de 2 à k :

$$T^{i}(N) = F_n^{i}(N + A_n) - A_n$$

Dans ce qui suit , on montre que la propriété est vraie pour $i = k + 1$. On écrit :

$$F_n^{k+1}(N + A_n) - A_n = F_n^{1}\left(F_n^{k}(N + A_n)\right) - A_n$$

Remplaçons $F_n^{k}(N + A_n)$ par $T^{k}(N) + A_n$ , on obtient :

$$F_n^{k+1}(N + A_n) - A_n = F_n^{1}\left(T^{k}(N) + A_n\right) - A_n$$

En utilisant la relation suivante

$$F_n^{1}(N + A_n) = T^{1}(N) + A_n$$

Nous permet d'écrire :

$$F_n^{1}\left(T^{k}(N) + A_n\right) = T^{1}\left(T^{k}(N)\right) + A_n$$

Par conséquent :

$$F_n^{1}\left(T^{k}(N) + A_n\right) - A_n = T^{1}\left(T^{k}(N)\right) = T^{k+1}(N)$$

Comme on a :

$$F_n^{1}\left(T^{k}(N) + A_n\right) - A_n = F_n^{1}\left(T^{k}(N) + A_n\right)$$

On déduit que :

$$F_n^{k+1}(N + A_n) - A_n = T^{k+1}(N)$$

Par suite la propriété est vraie pour i=k+1. On peut conclure que la propriété est vraie pour tout entier naturel non nul k

$$T^{k}(N) = F_n^{k}(N + A_n) - A_n$$

Pour obtenir l'autre équation il suffit de poser que :

$$P = N + A_n$$

**Proposition 2.2**

Soit n un entier relatif alors pour tous entiers naturels non nuls N et k, on a :

(2.8) $$T^{k}(N) = \frac{F_n^{k}(N + A_n) + F_{-n-1}^{k}(N + A_{-n-1})}{2}$$

**Démonstration :**

On sait que :

$$\begin{cases} T^{k}(N) = F_n^{k}(N + A_n) - A_n \\ T^{k}(N) = F_{-n-1}^{k}(N + A_{-n-1}) - A_{-n-1} \end{cases}$$

Par conséquent

$$2T^k(N) = F_n^k(N + A_n) - A_n + F_{-n-1}^k(N + A_{-n-1}) - A_{-n-1}$$
$$= F_n^k(N + A_n) + F_{-n-1}^k(N + A_{-n-1}) - (A_n + A_{-n-1})$$

Comme :

$$A_{-n-1} = 2(-n - 1) + 1$$
$$= -2n - 1$$

et on a aussi :

$$A_n = 2n + 1$$

Donc on peut déduire que :

$$A_{-n-1} + A_n = 0$$

Ceci nous peut de conclure que :

$$T^k(N) = \frac{F_n^k(N + A_n) + F_{-n-1}^k(N + A_{-n-1}) - (A_n + A_{-n-1})}{2}$$
$$= \frac{F_n^k(N + A_n) + F_{-n-1}^k(N + A_{-n-1})}{2}$$

**Théorème 2.2**

Soient n, k et N des entiers naturels non nuls alors on peur écrire :

(2.9) $$T^k(N) = \lim_{n\to+\infty}(\frac{1}{2(n+1)}\left(\sum_{j=0}^{n}((F_j^k\left(N + A_j\right) + F_{-j-1}^k\left(N + A_{-j-1}\right))\right))$$

Ce qui signifie que le terme $T^k(N)$ peut être interprété comme la moyenne arithmétique de tous les termes $(F_j^k(N + (2j + 1)))_{j\in\mathbb{Z}}$

**Démonstration**

D'après la proposition 2.3 on a :

$$T^k(N) = \frac{F_j^k(N + A_n) + F_{-j-1}^k\left(N + A_{-j-1}\right)}{2}$$

Qu'on peut l'écrire sous la forme suivante :

$$2T^k(N) = F_j^k\left(N + A_j\right) + F_{-j-1}^k\left(N + A_{-j-1}\right)$$

Cette égalité ne dépend pas de la valeur de l'entier j, on fait varier j de 0 à n (n un entier naturel non nul) comme suit :

$$\begin{cases} F_0^k(N+A_0) + F_{-1}^k(N+A_{-1}) = 2T^k(N) \\ F_1^k(N+A_1) + F_{-2}^k(N+A_{-2}) = 2T^k(N) \\ F_2^k(N+A_2) + F_{-3}^k(N+A_{-3}) = 2T^k(N) \\ \quad . \\ \quad . \\ F_{n-1}^k(N+A_{n-1}) + F_{-n}^k(N+A_{-n}) = 2T^k(N) \\ F_n^k(N+A_n) + F_{-n-1}^k(N+A_{-n-1}) = 2T^k(N) \end{cases}$$

Ce système qui contient (n+1) équations est équivalent à l'équation suivante :

$$2\,(n+1)T^k(N) = \sum_{j=0}^{n} (F_j^k(N+A_j) + F_{-j-1}^k(N+A_{-j-1}))$$

Ou encore :

$$T^k(N) = \frac{1}{2\,(n+1)} (\sum_{j=0}^{n} \left(F_j^k(N+A_j) + F_{-j-1}^k(N+A_{-j-1})\right))$$

Comme cette équation est valable quelque soit n donc on peut écrire :

$$T^k(N) = \lim_{n\to+\infty} (\frac{1}{2\,(n+1)} (\sum_{j=0}^{n} (F_j^k(N+A_j) + F_{-j-1}^k(N+A_{-j-1}))))$$

**Exemple 2.1**

Sur la figure suivante, on représente la variation de $F_n^k(N+A_n)$ et de $T^k(N)$ dans le cas ou N=9 et n=-3,n=-2,n=-1, n=0, n=1 et n=2.

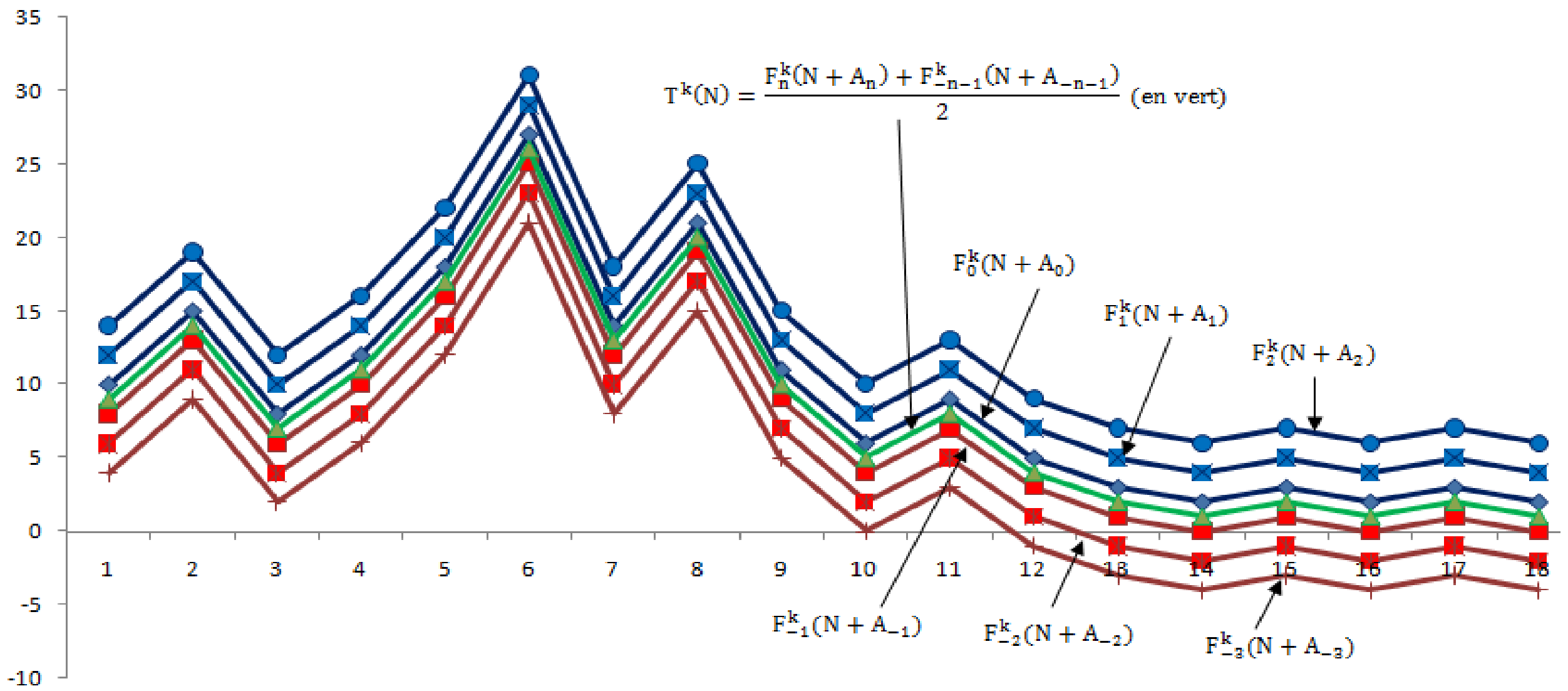

Figure 5: Comportement de différentes suites $S_{F,n}(N+A_n)$ et de la suite de Collatz Sy(N) pour N=9

On remarque bien que la suite de Collatz prend les valeurs moyennes des suites $S_{F,n}(N + A_n)$ et $S_{F,-n-1}(N + A_{-n-1})$

**Remarque 2.2**

Cette remarque est basée sur tous les résultats obtenus dans cet article qu'on peut le résumer dans le tableau suivant :

Tableau1: Illustration de la relation entre la fonction T et les deux fonctions $F_j$ et $F_{-j-1}$.

| | $F_{-j-1}$ | $F_j$ |
|---|---|---|
| n | $-j-1$ | j |
| Moyenne des indices | $\frac{-j-1+j}{2} = -\frac{1}{2}$ | |
| $A_n = 2n + 1$ | $-2j-1$ | $2j+1$ |
| Moyenne de $A_{-j-1}$ et $A_j$ | $\frac{A_{-j-1} + A_j}{2} = 0 = A_{-\frac{1}{2}}$ | |
| $C_{F,n} = (2n + 2{,}2n + 3)$ | $(-2j, -2j + 1)$ | $(2j + 2{,}2j + 3)$ |
| Moyenne de deux cycles $C_{F,-j-1}$ et $C_{F,j}$ | $\left(\frac{-2j + (2j+2)}{2}; \frac{(-2j+1) + (2j+3)}{2}\right) = (1{,}2) = C_{T,-\frac{1}{2}}$ | |
| Moyenne de deux termes $F^k_{-j-1}$ et $F^k_{-j}$ | $\frac{F^k_{-j-1}(N + A_{-j-1}) + F^k_{-j}(N + A_j)}{2} = T^k(N)$ | |

La figure suivante illustre cette propriété particulière de la fonction de Collatz qui apparait comme une fonction moyenne de toutes les fonctions de la suite considérée.

| | | $F_{-4}$ | | $F_{-3}$ | | $F_{-2}$ | | $F_{-1}$ | | T | | $F_0$ | | $F_1$ | | $F_2$ | | $F_3$ | | |
|---|---|---|---|---|---|---|---|---|---|---|---|---|---|---|---|---|---|---|---|---|
| | | | | | | | | | | | | | | | | | | | | |
| | | | | | | | | | | | | | | | | | | | | |
| 3 | | | | | | | | | | | | | | | | | | 8 | 9 | |
| 2 | | | | | | | | | | | | | | | | 6 | 7 | | | |
| 1 | | | | | | | | | | | | | | 4 | 5 | | | | | |
| 0 | | | | | | | | | | | | 2 | 3 | | | | | | | |
| $-\frac{1}{2}$ | | | | | | | | | | 1 | 2 | | | | | | | | | |
| -1 | | | | | | | | 0 | 1 | | | | | | | | | | | |
| -2 | | | | | | -2 | -1 | | | | | | | | | | | | | |
| -3 | | | | -4 | -3 | | | | | | | | | | | | | | | |
| -4 | | -6 | -5 | | | | | | | | | | | | | | | | | |
| | | | | | | | | | | | | | | | | | | | | |
| | | | | | | | | | | | | | | | | | | | | |

Figure 6 : Variation des domaines de validation de la conjecture de Collatz (en bleu) selon la fonction $F_n$ et exceptionnalité de la fonction T de Collatz

La fonction de Collatz T peut être interprétée comme la moyenne arithmétique de toutes les fonctions de la suite des fonctions $(F_n)_{n \in \mathbb{Z}}$

**Corollaire 2.1**

Soient N et k deux entiers naturels non nuls alors $\forall\, n \in \mathbb{Z}$

(2.10) $$F_n^k(N + A_n) - A_n = \text{constante}$$

**Exemple 2.2**

Prenons le cas de N=15 et k=3 sur le tableau suivant, on représente les différentes valeurs de $F_n^3(15 + A_n) - A_n$ pour quelques valeurs de n :

Tableau2 : Vérification de l'équation (2.10) pour quelques valeurs de n

| n | -4 | -3 | -2 | -1 | 0 | 1 | 2 | 3 | 4 |
|---|---|---|---|---|---|---|---|---|---|
| $A_n$ | -7 | -5 | -3 | -1 | 1 | 3 | 5 | 7 | 9 |
| $F_n^3(15 + A_n) - A_n$ | 53 | 53 | 53 | 53 | 53 | 53 | 53 | 53 | 53 |

On remarque que $F_n^3(15 + A_n) - A_n = T^3(15)$

**Corollaire 2.2**

Soient m et n deux entiers relatifs bien déterminés alors $\forall\, k \in \mathbb{N}^*$ et $\forall\, N \in \mathbb{N}^*$ on a :

(2.11) $$F_n^k(N + A_n) - F_m^k(N + A_m) = A_n - A_m = \text{constante}$$

Autrement la différence $F_n^k(N + A_n) - F_m^k(N + A_m)$ ne dépend pas ni de N ni de k, elle dépend uniquement de m et n.

**Démonstration**

D'après le théorème précédent (2.1), on peut écrire :

$$\begin{cases} F_n^k(N + A_n) - A_n = T^k(N) \\ F_m^k(N + A_m) - A_m = T^k(N) \end{cases}$$

Ce qui implique que :

$$F_n^k(N + A_n) - A_n = F_m^k(N + A_m) - A_m$$

Il en résulte que :

$$F_n^k(N + A_n) - F_m^k(N + A_m) = A_n - A_m$$

Pour deux valeurs bien déterminées de n et m la différence $(A_n - A_m)$ est une constante et par suite la différence $F_n^k(N + A_n) - F_m^k(N + A_m)$ ne dépend pas de k et de N , elle dépend de m et n.

**Exemple 2.3**

On prend le cas de $n = 3$ et $m = -2$ par exemple, ils correspondent aux deux fonctions suivantes :

$$F_3(P) = \begin{cases} \frac{3}{2}P - 3 & \text{si } P \equiv 0 \text{ (mod2)} \\ \frac{1}{2}P + \frac{7}{2} & \text{si } P \equiv 1 \text{ (mod2)} \end{cases} \quad ; \; F_{-2}(P) = \begin{cases} \frac{3}{2}P + 2 & \text{si } P \equiv 0 \text{ (mod2)} \\ \frac{1}{2}P - \frac{3}{2} & \text{si } P \equiv 1 \text{ (mod2)} \end{cases}$$

Les valeurs de $A_3$ et de $A_{-2}$ sont comme suit :

$$\begin{cases} A_3 = 2x3 + 1 = 7 \\ A_{-2} = -2x2 + 1 = -3 \end{cases}$$

On déduit la valeur de $A_3 - A_{-2}$ comme suit :

$$A_3 - A_{-2} = 10$$

D'après le corollaire (2.2) on peut écrire :

$$F_3^k(N + 7) - F_{-2}^k(N - 3) = 10$$

Prenons N=27 par exemple, dans ce cas on a:

$$N + A_3 = 27 + 7 = 34 \text{ et } N + A_{-2} = 27 - 3 = 24$$

On construit un tableau de dimensions 3x11 : 3 lignes et 11 colonnes tel que :

-La première ligne correspond à la suite $S_{F,3}(34,11)$.

-La deuxième ligne correspond à la suite $S_{F,-2}(24,11)$.

-La troisième ligne correspond à la différence de deux termes $F_3^k(27 + A_3) - F_{-2}^k(27 + A_1)$

Les résultats obtenus sont les suivants :

Tableau3 : Vérification de l'équation (2.11) pour N=27, n=3 et m=-2

| | | | | | | | | | | |
|---|---|---|---|---|---|---|---|---|---|---|
| $S_{F,3}(34,11)$ | 34 | 48 | 69 | 38 | 54 | 78 | 114 | 168 | 249 | 128 |
| $S_{F,-2}(24,11)$ | 24 | 38 | 59 | 28 | 44 | 68 | 104 | 158 | 239 | 118 |
| $F_3^k(34) - F_{-2}^k(24)$ | 10 | 10 | 10 | 10 | 10 | 10 | 10 | 10 | 10 | 10 |

On remarque bien que $F_3^k(27 + A_3) - F_{-2}^k(27 + A_{-2}) = A_3 + A_{-2} = 10$ quelque soit la valeur de k.

Un autre exemple pour N=22, dans ce cas on a:

$$N + A_3 = 22 + 7 = 29 \text{ et } N + A_{-2} = 22 - 3 = 19$$

Les résultats obtenus pour ce cas sont représentés dans le tableau suivant :

Tableau4 : Vérification de l'équation (2.11) pour N=22, n=3 et m=-2

| | | | | | | | | | | |
|---|---|---|---|---|---|---|---|---|---|---|
| $S_{F,3}(29,11)$ | 18 | 24 | 33 | 20 | 27 | 17 | 12 | 15 | 11 | 9 |
| $S_{F,-2}(19,11)$ | 8 | 14 | 23 | 10 | 17 | 7 | 2 | 5 | 1 | -1 |
| $F_3^k(27) - F_{-2}^k(17)$ | 10 | 10 | 10 | 10 | 10 | 10 | 10 | 10 | 10 | 10 |

**Corollaire 2.3**

Soient N un entier naturels non nul et n et m des entiers relatifs quelconques et k un entier naturel. On considère les deux suites $S_{F,n}(N + A_n, k)$ et $S_{F,m}(N + A_m, k)$ alors :

(2.12) $$\text{si } \lim_{k\to+\infty} F_n^k(N + A_n) = +\infty \Leftrightarrow \lim_{k\to+\infty} F_m^k(N + A_m) = +\infty$$

(2.13) $$\text{si } \lim_{k\to+\infty} F_n^k(N + A_n) = -\infty \Leftrightarrow \lim_{k\to+\infty} F_m^k(N + A_m) = -\infty$$

Si la suite $(F_n^k(N + A_n))_{k\geq 0}$ est bornée $\Leftrightarrow$ la suite $(F_m^k(N + A_m))_{k\geq 0}$ est bornée

**Démonstration**

C'est une conséquence directe de l'équation suivante déduite de l'équation 2.9 de corollaire (2.2):

$$F_n^k(N + A_n) = A_n - A_m + F_m^k(N + A_m)$$
$$= 2(n + m) + 2 + F_m^k(N + A_m)$$

**Théorème 2.3**

On suppose que la conjecture de Collatz est vérifiée pour tout entier naturel positif N c'est à dire que :

$\forall$ N $\in \mathbb{N}^*$, il existe un entier naturel non nul k tel que $T^k(N) = 1$

Alors pour tout entier relatif P tel que $\mathbf{P} \geq (\mathbf{2n} + \mathbf{2})$ avec n un entier relatif quelconque, il existe un entier naturel non nul m tel que :

(2.14) $$F_n^m(P) = 2n + 2$$

Inversement, on suppose que pour tout entier relatif P tel que $P \geq 2n + 2$, il existe un entier naturel non nul m tel que :

$$F_n^m(P) = 2n + 2$$

Alors $\forall$ N $\in \mathbb{N}^*$, il existe un entier naturel non nul k tel que $T^k(N) = 1$

**Démonstration**

On suppose que la conjecture de Collatz est vraie c'est à dire qu'elle est vérifiée pour tout entier naturel non nul N alors $\forall$ N $\in \mathbb{N}^*$, il existe un entier naturel non nul k tel que $T^k(N) = 1$

Comme on a :

$$T^k(N) = F_n^k(N + A_n) - A_n$$

On déduit que

$$F_n^k(N + A_n) - A_n = 1$$

Equivalente à :

$$F_n^k(N + A_n) = 1 + A_n$$

Remplaçons $A_n$ par 2n+1 , on peut conclure que :

$\forall$ N $\in \mathbb{N}^*$, il existe un entier naturel non nul k tel que $F_n^k(N + 2n + 1) = 2n + 2$

Posons :

$$N + 2n + 1 = P$$

Comme $N \geq 1$ alors $P \geq 2n + 2$

On peut conclure par suite que pour tout entier relatif P tel que $P \geq 2n + 2$, il existe un entier naturel non nul k tel que :

$$F_n^k(P) = 2n + 2$$

Inversement, on suppose que pour tout entier relatif P tel que $P \geq 2n + 2$, il existe un entier naturel non nul k tel que :

$$F_n^k(P) = 2n + 2$$

Ce qui nous permet d'écrire :

$$T^k(P - A_n) + A_n = F_n^k(P) = 2n + 2$$

Posons $N = P - A_n$ alors puisque $P \geq 2n + 2$ et $A_n = 2n + 1$ donc $N \geq 1$ ainsi on peut conclure que :

$\forall$ N $\in \mathbb{N}^*$, il existe un entier naturel non nul k tel que $T^k(N) = F_n^k(P) - A_n = 1$

**Corollaire 2.4**

Pour tout entier relatif n, les deux termes $T^k(N)$ et $F_n^k(N + A_n)$ ont des parités différentes autrement :

Si $T^k(N)$ est pair alors $F_n^k(N + A_n)$ est impair et si $T^k(N)$ est impair alors $F_n^k(N + A_n)$ est pair.

**Démonstration**

Ce corollaire est une conséquence de l'équation suivante :

$$T^k(N) = F_n^k(N + A_n) - A_n$$

Puisque $A_n$ est impair donc nécessairement $T^k(N)$ et $F_n^k(N + A_n)$ n'ont pas la même parité.

**Corollaire 2.5**

Soient N et k deux entiers naturels non nuls et n un entier relatif alors le nombre des entiers **impairs** dans une suite $Sy(N, k)$ qui est noté $\beta_k(N)$ égal au nombre des entiers **pairs** dans la suite $S_{F,n}(N + A_n, k)$. Cette propriété se traduit par l'équation suivante :

(2.15) $$\beta_k(N) = \alpha_{n,k}(N + A_n)$$

**Démonstration**

C'est une conséquence directe de corollaire (2.4) en effet si $T^k(N)$ est impair $i_k(N) = 1$ et comme dans ce cas $F_n^k(N + A_n)$ est pair alors $f_{n,k}(N) = 1$, on peut écrire alors :

$$i_k(N) = 1 \Leftrightarrow f_{n,k}(N) = 1$$

De même pour l'autre cas si $T^k(N)$ est pair donc $i_k(N) = 0$ et $F_n^k(N + A_n)$ est impair $f_{n,k}(N) = 0$, ce qui nous permet d'écrire :

$$i_k(N) = 0 \Leftrightarrow f_{n,k}(N) = 0$$

On déduit que :

$$\beta_k(N) = \sum_{j=0}^{k-1} i_j(N) = \sum_{j=0}^{k-1} f_{n,j}(N) = \alpha_{n,k}(N)$$

**Corollaire 2.6**

Soient n ,m des entiers relatifs, N et k deux entiers naturels non nuls alors on a :

(2.16) $$\frac{3^{\beta_k(N)}}{2^k}(A_n - A_m) + \varphi_{n,k}(N + A_n) - \varphi_{m,k}(N + A_m) = (A_n - A_m)$$

**Démonstration**

On sait que :

$$\begin{cases} F_m^k(N + A_m) = \dfrac{3^{\alpha_{m,k}(N+A_m)}}{2^k}(N + A_m) + \varphi_{m,k}(N + A_m) \\ F_n^k(N + A_n) = \dfrac{3^{\alpha_{n,k}(N+A_n)}}{2^k}(N + A_n) + \varphi_{n,k}(N + A_n) \end{cases}$$

Or

$$\alpha_{n,k}(N+A_n) = \alpha_{m,k}(N+A_m) = \beta_k(N)$$

$$T^k(N) = F_m^k(N+A_m) - A_m = F_n^k(N+A_n) - A_n$$

Alors on déduit que :

$$\frac{3^{\beta_k(N)}}{2^k}(N+A_n) + \varphi_{n,k}(N+A_n) - A_n = \frac{3^{\beta_k(N)}}{2^k}(N+A_m) + \varphi_{m,k}(N+A_m) - A_m$$

D'ou

$$\frac{3^{\beta_k(N)}}{2^k}(A_n - A_m) + \varphi_{n,k}(N+A_n) - \varphi_{m,k}(N+A_m) = (A_n - A_m)$$

**Corollaire 2.7**

Soit N un entier naturel non nul , n et m deux entiers relatifs quelconques donc on peut écrire:

(2.17) $$\text{si } \lim_{k\to+\infty} \frac{3^{\beta_k(N)}}{2^k} = 0 \Leftrightarrow \lim_{k\to+\infty} (\varphi_{n,k}(N+A_n) - \varphi_{m,k}(N+A_m)) = 2(n-m)$$

On suppose que ($n > m$) alors :

(2.18) $$\text{si } \lim_{k\to+\infty} \frac{3^{\beta_k(N)}}{2^k} = +\infty \Leftrightarrow \lim_{k\to+\infty} \left(\varphi_{n,k}(N+A_n) - \varphi_{m,k}(N+A_m)\right) = -\infty$$

**Proposition 2.3**

Soient n un entier relatif, N et k deux entiers naturels non nuls et on considère les équations suivantes :

(2.19) $$\begin{cases} T^k(N) = \dfrac{3^{\beta_k(N)}}{2^k} N + r_k(N) \\ F_n^k(N+A_n) = \dfrac{3^{\alpha_{n,k}(N+A_n)}}{2^k}(N+A_n) + \varphi_{n,k}(N+A_n) \end{cases}$$

Donc on peut admettre les relations suivantes

(2.20) $$\begin{cases} \dfrac{3^{\beta_k(N)}}{2^k} = \dfrac{3^{\alpha_{n,k}(N+A_n)}}{2^k} \\ r_k(N) = \varphi_{n,k}(N+A_n) + A_n \left(\dfrac{3^{\beta_k(N)}}{2^k} - 1\right) \end{cases}$$

**Démonstration**

Les expressions de $T^k(N)$ en fonction de N et k et l'expression de $F_n^k(N+A_n)$ en fonction de N et de k sont les suivantes :

$$T^k(N) = \frac{3^{\beta_k(N)}}{2^k}N + r_k(N)$$

$$F_n^k(N + A_n) = \frac{3^{\alpha_{n,k}(N+A_n)}}{2^k}(N + A_n) + \varphi_{n,k}(N + A_n)$$

D'après le premier théorème (2.1):

$$T^k(N) = F_n^k(N + A_n) - A_n$$

Remplaçons chaque terme par son expression, on obtient :

$$\frac{3^{\beta_k(N)}}{2^k}N + r_k(N) = \frac{3^{\alpha_{n,k}(N+A_n)}}{2^k}(N + A_n) + \varphi_{n,k}(N + A_n) - A_n$$

Comme on a :

$$\frac{3^{\beta_k(N)}}{2^k} = \frac{3^{\alpha_{n,k}(N+A_n)}}{2^k}$$

Donc on peut déduire que:

$$r_k(N) = \varphi_{n,k}(N + A_n) + A_n(\frac{3^{\beta_k(N)}}{2^k} - 1)$$

**Corollaire 2.8**

Soit n un entier relatif quelconque donc $\forall$ P $\in \mathbb{Z}$ tel que $P \geq 2n + 2$ et $\forall k \in \mathbb{N}$ on peut écrire :

(2.21) $$F_n^k(P \geq 2n + 2) \geq 2n + 2$$

Autrement les suites $(F_n^k(P))_{k\geq 0, P\geq 2n+2}$ pour tout entier relatif $P \geq 2n + 2$ sont minorées par (2n+2)

**Démonstration**

On sait que $\forall$ N $\in \mathbb{N}^*$ et $\forall k \in \mathbb{N}$ on peut écrire :

$$T^k(N) \geq 1$$

Il en résulte que :

$$F_n^k(N + A_n) - A_n \geq 1$$

Ou encore :

$$F_n^k(N + A_n) \geq 1 + A_n = 2n + 2$$

On a :

$$N \geq 1 \Rightarrow N + A_n \geq 2n + 2$$

On pose $P = N + A_n$

On peut conclure que pour tout entier relatif P tel que $P \geq 2n+2$ et $\forall\, k \in \mathbb{N}$ on peut écrire :

$$F_n^k(P) \geq 2n+2$$

## 3. La matrice généralisée de Collatz : $\mathbb{T}_{\mathbf{n}}(\mathbf{2n+2, 2n+M, k})$

**Définition 3.1**

Un ensemble homogène continu des suites est un ensemble des suites issues d'une même fonction, elles ont la même longueur et de plus leurs premiers termes décrivent une suite arithmétique de raison arithmétique 1. Par exemple, on considère l'ensemble suivant :

(3.1) $\mathbb{S}_{F,n}(2n+2, P, k) = (S_{F,n}(2n+2, k), S_{F,n}(2n+3, k), S_{F,n}(2n+4, k), \ldots, S_{F,n}(2n+P, k))$

Cet ensemble constitue un ECH des suites générées par la fonction $F_n$ , elles ont la même longueur k de plus leurs premiers termes constituent une suite arithmétique de rasion1.

$$D(n) = (2n+2, 2n+3, 2n+4, \ldots., 2n+P, \ldots)$$

### 3.1 Représentation matricielle d'une distribution homogène continue :

La matrice généralisée de Collatz est d'une importance particulière, elle se comporte comme un tableau de dimension finie ou infinie. La **première colonne** de ce tableau contient toujours une suite arithmétique de raison 1 tel que les termes de la suite contenue dans la première colonne sont exprimés en fonction de n comme suit

(3.2) $$(2n+2, 2n+3, 2n+4, \ldots, 2n+M)$$

Avec M un entier relatif tel que $M > (2n+2)$

Chaque ligne contient les termes de la suite générée par la fonction $F_n$ et dont le premier terme est l'entier placé dans la première case de cette ligne. Comme montre l'exemple suivant :

| | | | | | | |
|---|---|---|---|---|---|---|
| $2n+P$ | $F_n^1(2n+P)$ | $F_n^2(2n+P)$ | | | $F_n^{k-1}(2n+P)$ | $F_n^k(2n+P)$ |

La matrice généralisée de Collatz est notée $\mathbb{T}_{\mathbf{n}}(\mathbf{2n+2, 2n+M, k})$

Avec k la longueur des suites qui correspond au nombre des colonnes dans la matrice considérée. Sur la figure suivante, on représente la matrice généralisée de Collatz $\mathbb{T}_n(2n+2, 2n+16, 6)$

| $P_i = 2n + N_i$ | $F_n^1(P_i)$ | $F_n^2(P_i)$ | $F_n^3(P_i)$ | $F_n^4(P_i)$ | $F_n^5(P_i)$ |
|---|---|---|---|---|---|
| | | | | | |
| | | | | | |
| | | | | | |
| 2n + 16 | 2n + 24 | 2n + 36 | 2n + 54 | 2n + 81 | 2n + 41 |
| 2n + 15 | 2n + 8 | 2n + 12 | 2n + 18 | 2n + 27 | 2n + 14 |
| 2n + 14 | 2n + 21 | 2n + 11 | 2n + 6 | 2n + 9 | 2n + 5 |
| 2n + 13 | 2n + 7 | 2n + 4 | 2n + 6 | 2n + 9 | 2n + 5 |
| 2n + 12 | 2n + 18 | 2n + 27 | 2n + 14 | 2n + 21 | 2n + 11 |
| 2n + 11 | 2n + 6 | 2n + 9 | 2n + 5 | 2n + 3 | 2n + 2 |
| 2n + 10 | 2n + 15 | 2n + 8 | 2n + 12 | 2n + 18 | 2n + 27 |
| 2n + 9 | 2n + 5 | 2n + 3 | 2n + 2 | 2n + 3 | 2n + 2 |
| 2n + 8 | 2n + 12 | 2n + 18 | 2n + 27 | 2n + 14 | 2n + 21 |
| 2n + 7 | 2n + 4 | 2n + 6 | 2n + 9 | 2n + 5 | 2n + 3 |
| 2n + 6 | 2n + 9 | 2n + 5 | 2n + 3 | 2n + 2 | 2n + 3 |
| 2n + 5 | 2n + 3 | 2n + 2 | 2n + 3 | 2n + 2 | 2n + 3 |
| 2n + 4 | 2n + 6 | 2n + 9 | 2n + 5 | 2n + 3 | 2n + 2 |
| 2n + 3 | 2n + 2 | 2n + 3 | 2n + 2 | 2n + 3 | 2n + 2 |
| 2n + 2 | 2n + 3 | 2n + 2 | 2n + 3 | 2n + 2 | 2n + 3 |

Figure 7: Matrice généralisée de Collatz $\mathbb{T}_n(2n + 2,2n + 16,6)$

Pour Les valeurs de n, on peut distinguer deux cas :

$-$Si $n = -\frac{1}{2}$

Si on remplace n par $\left(-\frac{1}{2}\right)$ dans la matrice généralisée on obtient les suites generées par la fonction de collatz T. C'est la seule exception ou n peut prendre un entier non relatif.

Les suites contenues dans ce tableau sont des suites de Collatz. Dans ce cas, la matrice généralisée est comme suit :

(3.3) $$\mathbb{T}_n(2n + 2,2n + M, k) \to \mathbb{T}_{-\frac{1}{2}}(1, M - 1, k)$$

$-$Si $n \in \mathbb{Z}$

Dans ce cas n est un entier relatif bien déterminé et la matrice généralisée contient toutes les suites générées par la fonction $F_n$ et dont les premiers termes sont compris entre $2n + 2$ et $2n + M$.

## 3.2 Les règles de l'équivalence chromatique

Pour une matrice qui contient **des suites de Collatz** $(Sy(N,k))$ c'est à dire la matrice $\mathbb{T}_{-\frac{1}{2}}(1, P, k)$ les cases sont colorées selon les règles suivantes :

| |
|---|
| 4k |
| 4k + 1 |
| 4k + 2 |
| 4k + 3 |
| 4k |
| 4k + 1 |
| 4k + 2 |
| 4k + 3 |

Figure 8: Les règles de la coloration des cases selon la forme des entiers placés à l'intérieur de ces cases

Pour les autres matrices qui contiennent des suites $S_{F,n}(P, k)$ générées par une fonction $F_n$donc elles s'agissent des matrices $\mathbb{T}_{n\in\mathbb{Z}}(2n.2n + M, k)$ on parle dès règles d'équivalence chromatique entre les matrices $\mathbb{T}_{n\in\mathbb{Z}}(2n + 2, \dots)$ et la matrice $\mathbb{T}_{-\frac{1}{2}}(1, \dots)$

Les règles de l'équivalence chromatique consistent à colorer chaque case qui contient un terme sous forme $F_n^k(N + A_n)$ du tableau $\mathbb{T}_{n\in\mathbb{Z}}(2n + 2{,}2n + M, k)$ par la même couleur de la case qui contient le terme $T^k(N)$ dans le tableau $\mathbb{T}_{-\frac{1}{2}}(1, M - 1, k)$.

Par exemple si $T^k(N)$ s'écrit sous la forme 4j+1 alors la case qui contient ce terme sera colorée en bleu ce qui signifie que la case qui contient le terme $F_n^k(N + A_n)$ dans le tableau $\mathbb{T}_n(2n + 2{,}2n + M, k)$ sera colorée aussi par une couleur bleue.

Les règles de l'équivalence chromatique  sont déduites à partir de la matrice qui contient des suites de Collatz

| | | |
|---|---|---|
| 1 | 2 | 1 |
| 2 | 1 | 2 |
| 3 | 5 | 8 |
| 4 | 2 | 1 |
| 5 | 8 | 4 |
| 6 | 3 | 5 |
| 2n + 2 | 2n + 3 | 2n + 2 |
| 2n + 3 | 2n + 2 | 2n + 3 |
| 2n + 4 | 2n + 6 | 2n + 9 |
| 2n + 5 | 2n + 3 | 2n + 2 |
| 2n + 6 | 2n + 9 | 2n + 5 |
| 2n + 7 | 2n + 4 | 2n + 6 |

Dans ce cas les cases d'une matrice généralisée sont colorées selon la règle suivante dite d'équivalence chromatique :

| |
|---|
| 2n + 1 |
| 2n + 2 |
| 2n + 3 |
| 2n + 4 |

Figure 9: Les règles de la coloration des cases des tableaux contenants les suites générées par une fonction $F_n$

Ceci peut être interprété comme suit : On fait attribuer une couleur bleu à tous les entiers de forme 2n+2, une couleur jaune à tous les entiers qui s'écrivent sous la forme 2n+3, on fait attribuer une couleur rouge à tous les entiers dont les formes 2n+4 et pour les entier qui s'écrivent sous la forme 2n+1, les cases seront colorées en vert.

## 3.3 Application à la matrice généralisée de Collatz

On fait appliquer ces règles (figure 9) pour la coloration de différentes cases de la matrice généralisée de Collatz, on obtient la matrice chromatisée suivante :

| | | | | | |
|---|---|---|---|---|---|
| 2n + 16 | 2n + 24 | 2n + 36 | 2n + 54 | 2n + 81 | 2n + 41 |
| 2n + 15 | 2n + 8 | 2n + 12 | 2n + 18 | 2n + 27 | 2n + 14 |
| 2n + 14 | 2n + 21 | 2n + 11 | 2n + 6 | 2n + 9 | 2n + 5 |
| 2n + 13 | 2n + 7 | 2n + 4 | 2n + 6 | 2n + 9 | 2n + 5 |
| 2n + 12 | 2n + 18 | 2n + 27 | 2n + 14 | 2n + 21 | 2n + 11 |
| 2n + 11 | 2n + 6 | 2n + 9 | 2n + 5 | 2n + 3 | 2n + 2 |
| 2n + 10 | 2n + 15 | 2n + 8 | 2n + 12 | 2n + 18 | 2n + 27 |
| 2n + 9 | 2n + 5 | 2n + 3 | 2n + 2 | 2n + 3 | 2n + 2 |
| 2n + 8 | 2n + 12 | 2n + 18 | 2n + 27 | 2n + 14 | 2n + 21 |
| 2n + 7 | 2n + 4 | 2n + 6 | 2n + 9 | 2n + 5 | 2n + 3 |
| 2n + 6 | 2n + 9 | 2n + 5 | 2n + 3 | 2n + 2 | 2n + 3 |
| 2n + 5 | 2n + 3 | 2n + 2 | 2n + 3 | 2n + 2 | 2n + 3 |
| 2n + 4 | 2n + 6 | 2n + 9 | 2n + 5 | 2n + 3 | 2n + 2 |
| 2n + 3 | 2n + 2 | 2n + 3 | 2n + 2 | 2n + 3 | 2n + 2 |
| 2n + 2 | 2n + 3 | 2n + 2 | 2n + 3 | 2n + 2 | 2n + 3 |

Figure 10: La matrice généralisée de Collatz $\mathbb{T}_n(2n+2{,}2n+16{,}6)$

On fait appliquer ces règles à la matrice de base de Collatz $\mathbb{T}_{-\frac{1}{2}}(1{,}15{,}6)$, on obtient le tableau chromatisé suivant :

| | | | | | |
|---|---|---|---|---|---|
| 15 | 23 | 35 | 53 | 80 | 20 |
| 14 | 7 | 11 | 17 | 26 | 13 |
| 13 | 20 | 10 | 5 | 8 | 4 |
| 12 | 6 | 3 | 5 | 8 | 4 |
| 11 | 17 | 26 | 13 | 20 | 10 |
| 10 | 5 | 8 | 4 | 2 | 1 |
| 9 | 14 | 7 | 11 | 17 | 26 |
| 8 | 4 | 2 | 1 | 2 | 1 |
| 7 | 11 | 17 | 26 | 13 | 20 |
| 6 | 3 | 5 | 8 | 4 | 2 |
| 5 | 8 | 4 | 2 | 1 | 2 |
| 4 | 2 | 1 | 2 | 1 | 2 |
| 3 | 5 | 8 | 4 | 2 | 1 |
| 2 | 1 | 2 | 1 | 2 | 1 |
| 1 | 2 | 1 | 2 | 1 | 1 |

Figure 11: La matrice de base de Collatz $\mathbb{T}_{-\frac{1}{2}}(1{,}15{,}6)$

Dans ce qui suit on représente les deux matrices $\mathbb{T}_{-\frac{1}{2}}(1,15,6)$ et $\mathbb{T}_{n}(2n+2,2n+16,6)$ sur la même figure.

| | | | | | |
|---|---|---|---|---|---|
| 15 | 23 | 35 | 53 | 80 | 20 |
| 14 | 7 | 11 | 17 | 26 | 13 |
| 13 | 20 | 10 | 5 | 8 | 4 |
| 12 | 6 | 3 | 5 | 8 | 4 |
| 11 | 17 | 26 | 13 | 20 | 10 |
| 10 | 5 | 8 | 4 | 2 | 1 |
| 9 | 14 | 7 | 11 | 17 | 26 |
| 8 | 4 | 2 | 1 | 2 | 1 |
| 7 | 11 | 17 | 26 | 13 | 20 |
| 6 | 3 | 5 | 8 | 4 | 2 |
| 5 | 8 | 4 | 2 | 1 | 2 |
| 4 | 2 | 1 | 2 | 1 | 2 |
| 3 | 5 | 8 | 4 | 2 | 1 |
| 2 | 1 | 2 | 1 | 2 | 1 |
| 1 | 2 | 1 | 2 | 1 | 1 |
| 2n + 16 | 2n + 24 | 2n + 36 | 2n + 54 | 2n + 81 | 2n + 41 |
| 2n + 15 | 2n + 8 | 2n + 12 | 2n + 18 | 2n + 27 | 2n + 14 |
| 2n + 14 | 2n + 21 | 2n + 11 | 2n + 6 | 2n + 9 | 2n + 5 |
| 2n + 13 | 2n + 7 | 2n + 4 | 2n + 6 | 2n + 9 | 2n + 5 |
| 2n + 12 | 2n + 18 | 2n + 27 | 2n + 14 | 2n + 21 | 2n + 11 |
| 2n + 11 | 2n + 6 | 2n + 9 | 2n + 5 | 2n + 3 | 2n + 2 |
| 2n + 10 | 2n + 15 | 2n + 8 | 2n + 12 | 2n + 18 | 2n + 27 |
| 2n + 9 | 2n + 5 | 2n + 3 | 2n + 2 | 2n + 3 | 2n + 2 |
| 2n + 8 | 2n + 12 | 2n + 18 | 2n + 27 | 2n + 14 | 2n + 21 |
| 2n + 7 | 2n + 4 | 2n + 6 | 2n + 9 | 2n + 5 | 2n + 3 |
| 2n + 6 | 2n + 9 | 2n + 5 | 2n + 3 | 2n + 2 | 2n + 3 |
| 2n + 5 | 2n + 3 | 2n + 2 | 2n + 3 | 2n + 2 | 2n + 3 |
| 2n + 4 | 2n + 6 | 2n + 9 | 2n + 5 | 2n + 3 | 2n + 2 |
| 2n + 3 | 2n + 2 | 2n + 3 | 2n + 2 | 2n + 3 | 2n + 2 |
| 2n + 2 | 2n + 3 | 2n + 2 | 2n + 3 | 2n + 2 | 2n + 3 |

Figure 12: Comportements des suites dans les deux matrices $\mathbb{T}_{-\frac{1}{2}}(1,15,6)$ et $\mathbb{T}_{n}(2n+2,2n+16,6)$

La principale interprétation est la suivante : les termes qui s'écrivent sous la forme 2n+4j dans la matrice généralisée possèdent la même distribution des termes qui s'écrivent sous la forme 4j+3 dans la matrice de Collatz, les termes qui prennent la forme (2n+4j+2) possèdent la distribution que les termes dont la forme est (4j+1) , les termes qui ont la forme (2n+4j+1) possèdent la même distribution que les termes dont la forme est 4j dans la matrice de Collatz alors que les termes qui s'écrivent sous la forme (2n+4j+3 )possèdent la même distribution que les termes qui s'écrivent sous la forme (4j+2) dans la matrice de Collatz.

On peut vérifier que la matrice $\mathbb{T}_{-\frac{1}{2}}(1,15,6)$ peut être déduite à partir de la matrice $\mathbb{T}_{n}(2n+2,2n+16,6)$ on fait remplacer n par $(-0.5)$.

## 4. Quelques exemples des suites générées par des fonctions $(\mathbf{F_n})_{n\in\mathbb{Z}}$

### 4.1 Comportements des suites $\mathbf{S_{F0}(P,n)}$

Les suites $S_{F,0}(P,n)$ sont générées par la fonction $F_0$ qui est définie sur D(0) comme suit :

(4.1) $$F_0(P) = \begin{cases} \frac{3}{2}P & \text{si } P \equiv 0 (\text{mod}2) \\ \frac{P}{2}+\frac{1}{2} & \text{si } P \equiv 1\ (\text{mod}2) \end{cases}$$

Les points particuliers ou caractéristiques de cette fonction sont les suivants (n=0) :

$$\begin{cases} P_{0,1} = 2 \\ P_{0,2} = 3 \\ A_0 = 1 \end{cases} \quad (4.2)$$

On peut vérifier que :

$$F_0(2) = 3 \,;\, F_0(3) = 2$$

Sur la figure suivante, on représente la matrice $\mathbb{T}_0(2{,}17{,}10)$ qui contient les 16 premières suites $S_{F,0}(P, 11)$ de même longueur 11. Les cases sont colorées selon les règles indiquées sur la figure 9 Pour obtenir cette matrice, il suffit de remplacer n par 0 dans la matrice généralisée de Collatz donnée par la figure 8.

| | | | | | | | | | | |
|---|---|---|---|---|---|---|---|---|---|---|
| 2 | 3 | 2 | 3 | 2 | 3 | 2 | 3 | 2 | 3 | 2 |
| 3 | 2 | 3 | 2 | 3 | 2 | 3 | 2 | 3 | 2 | 3 |
| 4 | 6 | 9 | 5 | 3 | 2 | 3 | 2 | 3 | 2 | 3 |
| 5 | 3 | 2 | 3 | 2 | 3 | 2 | 3 | 2 | 3 | 2 |
| 6 | 9 | 5 | 3 | 2 | 3 | 2 | 3 | 2 | 3 | 2 |
| 7 | 4 | 6 | 9 | 5 | 3 | 2 | 3 | 2 | 3 | 2 |
| 8 | 12 | 18 | 27 | 14 | 21 | 11 | 6 | 9 | 5 | 3 |
| 9 | 5 | 3 | 2 | 3 | 2 | 3 | 2 | 3 | 2 | 3 |
| 10 | 15 | 8 | 12 | 18 | 27 | 14 | 21 | 11 | 6 | 9 |
| 11 | 6 | 9 | 5 | 3 | 2 | 3 | 2 | 3 | 2 | 3 |
| 12 | 18 | 27 | 14 | 21 | 11 | 6 | 9 | 5 | 3 | 2 |
| 13 | 7 | 4 | 6 | 9 | 5 | 3 | 2 | 3 | 2 | 3 |
| 14 | 21 | 11 | 6 | 9 | 5 | 3 | 2 | 3 | 2 | 3 |
| 15 | 8 | 12 | 18 | 27 | 14 | 21 | 11 | 6 | 9 | 5 |
| 16 | 24 | 36 | 54 | 81 | 41 | 21 | 11 | 6 | 9 | 5 |
| 17 | 9 | 5 | 3 | 2 | 3 | 2 | 3 | 2 | 3 | 2 |

Figure 13 : La matrice $\mathbb{T}_0(2{,}17{,}11)$

On remarque bien que ces suites atteignent le cycle (2,3) après un nombre finie des itérations. On croit que cette propriété est vraie pour tout entier naturel N tel que $N \geq 2$, cette interprétation est basée sur le théorème 2.1 décrivant la relation entre $T^k(N)$ et $T^k(N + A_n)$ on parle dans ce cas la conjecture de Collatz pour la fonction $F_0$. Cette conjecture peut être énoncée comme suit :

$\forall$ P $\in \mathbb{N}$ tel que $P \geq 2$, il existe un entier naturel non nul k tel que :

$$F_0^k(P) = 2 \quad (4.3)$$

On fait comparer les deux matrices $\mathbb{T}_0(2{,}17{,}11)$ et $\mathbb{T}_{-\frac{1}{2}}(1, 16, 11)$ a fin de mettre en relief les propriétés communes et de similarité entre les comportements des différentes suites générées par les fonctions T et $F_0$. Les cases sont colorées tout en respectant les règles de l'équivalence chromatique de la figure 9.

| | | | | | | | | | | |
|---|---|---|---|---|---|---|---|---|---|---|
| 1 | 2 | 1 | 2 | 1 | 2 | 1 | 2 | 1 | 2 | 1 |
| 2 | 1 | 2 | 1 | 2 | 1 | 2 | 1 | 2 | 1 | 2 |
| 3 | 5 | 8 | 4 | 2 | 1 | 2 | 1 | 2 | 1 | 2 |
| 4 | 2 | 1 | 2 | 1 | 2 | 1 | 2 | 1 | 2 | 1 |
| 5 | 8 | 4 | 2 | 1 | 2 | 1 | 2 | 1 | 2 | 1 |
| 6 | 3 | 5 | 8 | 4 | 2 | 1 | 2 | 1 | 2 | 1 |
| 7 | 11 | 17 | 26 | 13 | 20 | 10 | 5 | 8 | 4 | 2 |
| 8 | 4 | 2 | 1 | 2 | 1 | 2 | 1 | 2 | 1 | 2 |
| 9 | 14 | 7 | 11 | 17 | 26 | 13 | 20 | 10 | 5 | 8 |
| 10 | 5 | 8 | 4 | 2 | 1 | 2 | 1 | 2 | 1 | 2 |
| 11 | 17 | 26 | 13 | 20 | 10 | 5 | 8 | 4 | 2 | 1 |
| 12 | 6 | 3 | 5 | 8 | 4 | 2 | 1 | 2 | 1 | 2 |
| 13 | 20 | 10 | 5 | 8 | 4 | 2 | 1 | 2 | 1 | 2 |
| 14 | 7 | 11 | 17 | 26 | 13 | 20 | 10 | 5 | 8 | 4 |
| 15 | 23 | 35 | 53 | 80 | 40 | 20 | 10 | 5 | 8 | 4 |
| 16 | 8 | 4 | 2 | 1 | 2 | 1 | 2 | 1 | 2 | 1 |
| 2 | 3 | 2 | 3 | 2 | 3 | 2 | 3 | 2 | 3 | 2 |
| 3 | 2 | 3 | 2 | 3 | 2 | 3 | 2 | 3 | 2 | 3 |
| 4 | 6 | 9 | 5 | 3 | 2 | 3 | 2 | 3 | 2 | 3 |
| 5 | 3 | 2 | 3 | 2 | 3 | 2 | 3 | 2 | 3 | 2 |
| 6 | 9 | 5 | 3 | 2 | 3 | 2 | 3 | 2 | 3 | 2 |
| 7 | 4 | 6 | 9 | 5 | 3 | 2 | 3 | 2 | 3 | 2 |
| 8 | 12 | 18 | 27 | 14 | 21 | 11 | 6 | 9 | 5 | 3 |
| 9 | 5 | 3 | 2 | 3 | 2 | 3 | 2 | 3 | 2 | 3 |
| 10 | 15 | 8 | 12 | 18 | 27 | 14 | 21 | 11 | 6 | 9 |
| 11 | 6 | 9 | 5 | 3 | 2 | 3 | 2 | 3 | 2 | 3 |
| 12 | 18 | 27 | 14 | 21 | 11 | 6 | 9 | 5 | 3 | 2 |
| 13 | 7 | 4 | 6 | 9 | 5 | 3 | 2 | 3 | 2 | 3 |
| 14 | 21 | 11 | 6 | 9 | 5 | 3 | 2 | 3 | 2 | 3 |
| 15 | 8 | 12 | 18 | 27 | 14 | 21 | 11 | 6 | 9 | 5 |
| 16 | 24 | 36 | 54 | 81 | 41 | 21 | 11 | 6 | 9 | 5 |
| 17 | 9 | 5 | 3 | 2 | 3 | 2 | 3 | 2 | 3 | 2 |

Figure 14: Comportements de différentes suites dans les deux matrice $\mathbb{T}_0(2,17,11)$ et $\mathbb{T}_{-\frac{1}{2}}(1,16,11)$

On remarque bien que la suite de premier terme un entier naturel N dans le tableau $\mathbb{T}_{-\frac{1}{2}}(1,16,10)$ possède la même distribution chromatique que la suite de premier terme $(N+1)$ dans la matrice $\mathbb{T}_0(2,17,11)$ on peut vérifier aisément la relation suivante :

$$F_0^k(N+1) - 1 = T^k(N) \tag{4.4}$$

## 4.2 Comportements des suites générées par la fonction $F_1$

La fonction $F_1$ (n=1) est définie sur D(1) comme suit :

$$F_1(P) = \begin{cases} \frac{3}{2}P - 1 & \text{si } P \equiv 0 (\text{mod} 2) \\ \frac{P}{2} + \frac{3}{2} & \text{si } P \equiv 1 \ (\text{mod} 2) \end{cases} \tag{4.5}$$

Les points particuliers de cette fonction sont :

$$\begin{cases} P_{1,1} = 4 \\ P_{1,2} = 5 \\ A_1 = 3 \end{cases} \tag{4.6}$$

La figure suivante correspond à la matrice $\mathbb{T}_1(4,19,10)$ elle est obtenue par remplacement de n par 1 dans la matrice généralisée tout en respectant les règles de la coloration selon le principe de l'équivalence chromatique sus-indiqué :

| | | | | | | | | | |
|---|---|---|---|---|---|---|---|---|---|
| 4 | 5 | 4 | 5 | 4 | 5 | 4 | 5 | 4 | 5 |
| 5 | 4 | 5 | 4 | 5 | 4 | 5 | 4 | 5 | 4 |
| 6 | 8 | 11 | 7 | 5 | 4 | 5 | 4 | 5 | 4 |
| 7 | 5 | 4 | 5 | 4 | 5 | 4 | 5 | 4 | 5 |
| 8 | 11 | 7 | 5 | 4 | 5 | 4 | 5 | 4 | 5 |
| 9 | 6 | 8 | 11 | 7 | 5 | 4 | 5 | 4 | 5 |
| 10 | 14 | 20 | 29 | 16 | 23 | 13 | 8 | 11 | 7 |
| 11 | 7 | 5 | 4 | 5 | 4 | 5 | 4 | 5 | 4 |
| 12 | 17 | 10 | 14 | 20 | 29 | 16 | 23 | 13 | 8 |
| 13 | 8 | 11 | 7 | 5 | 4 | 5 | 4 | 5 | 4 |
| 14 | 20 | 29 | 16 | 23 | 13 | 8 | 11 | 7 | 5 |
| 15 | 9 | 6 | 8 | 11 | 7 | 5 | 4 | 5 | 4 |
| 16 | 23 | 13 | 8 | 11 | 7 | 5 | 4 | 5 | 4 |
| 17 | 10 | 14 | 20 | 29 | 16 | 23 | 13 | 8 | 11 |
| 18 | 26 | 38 | 56 | 83 | 43 | 23 | 13 | 8 | 11 |
| 19 | 11 | 7 | 5 | 4 | 5 | 4 | 5 | 4 | 5 |

Figure 15: La matrice $\mathbb{T}_1(4,19,10)$

On peut vérifier que :

$$F_1(4) = 5;\ F_1(5) = 4$$

Sur la figure suivante, on représente les deux matrices équivalentes $\mathbb{T}_1(4,19,11)$ et $\mathbb{T}_{-\frac{1}{2}}(1,16,11)$

| | | | | | | | | | | |
|---|---|---|---|---|---|---|---|---|---|---|
| 4 | 5 | 4 | 5 | 4 | 5 | 4 | 5 | 4 | 5 | 4 |
| 5 | 4 | 5 | 4 | 5 | 4 | 5 | 4 | 5 | 4 | 5 |
| 6 | 8 | 11 | 7 | 5 | 4 | 5 | 4 | 5 | 4 | 5 |
| 7 | 5 | 4 | 5 | 4 | 5 | 4 | 5 | 4 | 5 | 4 |
| 8 | 11 | 7 | 5 | 4 | 5 | 4 | 5 | 4 | 5 | 4 |
| 9 | 6 | 8 | 11 | 7 | 5 | 4 | 5 | 4 | 5 | 4 |
| 10 | 14 | 20 | 29 | 16 | 23 | 13 | 8 | 11 | 7 | 5 |
| 11 | 7 | 5 | 4 | 5 | 4 | 5 | 4 | 5 | 4 | 5 |
| 12 | 17 | 10 | 14 | 20 | 29 | 16 | 23 | 13 | 8 | 11 |
| 13 | 8 | 11 | 7 | 5 | 4 | 5 | 4 | 5 | 4 | 5 |
| 14 | 20 | 29 | 16 | 23 | 13 | 8 | 11 | 7 | 5 | 4 |
| 15 | 9 | 6 | 8 | 11 | 7 | 5 | 4 | 5 | 4 | 5 |
| 16 | 23 | 13 | 8 | 11 | 7 | 5 | 4 | 5 | 4 | 5 |
| 17 | 10 | 14 | 20 | 29 | 16 | 23 | 13 | 8 | 11 | 7 |
| 18 | 26 | 38 | 56 | 83 | 43 | 23 | 13 | 8 | 11 | 7 |
| 19 | 11 | 7 | 5 | 4 | 5 | 4 | 5 | 4 | 5 | 4 |
| 1 | 2 | 1 | 2 | 1 | 2 | 1 | 2 | 1 | 2 | 1 |
| 2 | 1 | 2 | 1 | 2 | 1 | 2 | 1 | 2 | 1 | 2 |
| 3 | 5 | 8 | 4 | 2 | 1 | 2 | 1 | 2 | 1 | 2 |
| 4 | 2 | 1 | 2 | 1 | 2 | 1 | 2 | 1 | 2 | 1 |
| 5 | 8 | 4 | 2 | 1 | 2 | 1 | 2 | 1 | 2 | 1 |
| 6 | 3 | 5 | 8 | 4 | 2 | 1 | 2 | 1 | 2 | 1 |
| 7 | 11 | 17 | 26 | 13 | 20 | 10 | 5 | 8 | 4 | 2 |
| 8 | 4 | 2 | 1 | 2 | 1 | 2 | 1 | 2 | 1 | 2 |
| 9 | 14 | 7 | 11 | 17 | 26 | 13 | 20 | 10 | 5 | 8 |
| 10 | 5 | 8 | 4 | 2 | 1 | 2 | 1 | 2 | 1 | 2 |
| 11 | 17 | 26 | 13 | 20 | 10 | 5 | 8 | 4 | 2 | 1 |
| 12 | 6 | 3 | 5 | 8 | 4 | 2 | 1 | 2 | 1 | 2 |
| 13 | 20 | 10 | 5 | 8 | 4 | 2 | 1 | 2 | 1 | 2 |
| 14 | 7 | 11 | 17 | 26 | 13 | 20 | 10 | 5 | 8 | 4 |
| 15 | 23 | 35 | 53 | 80 | 40 | 20 | 10 | 5 | 8 | 4 |
| 16 | 8 | 4 | 2 | 1 | 2 | 1 | 2 | 1 | 2 | 1 |

Figure 16 : Les matrices $\mathbb{T}_1(4,19,11)$ et $\mathbb{T}_{-\frac{1}{2}}(1,16,11)$

On remarque que la relation suivante déjà démontrée dans le cas général est vérifiée pour tous les termes de deux tableaux considérés.

(4.7) $$F_1^k(N+3) = T^k(N) + 3$$

En se basant sur le théorème 2.1, la conjecture de Collatz peut être énoncée pour la fonction $F_1$ de la manière suivante :

$\forall$ P $\in \mathbb{N}$ tel que P $\geq$ 4, il existe un entier naturel non nul k tel que :

(4.8) $$F_1^k(P) = 4$$

En réalité on peut vérifier aisément que ces deux tableaux vérifiées toutes les propriétés déjà démontrées

## 4.3 Comportements des suites générées par la fonction $F_2$

Cette fonction correspond a une valeur n=2, elle est définie sur D(2) comme suit :

(4.9) $$F_2(P) = \begin{cases} \frac{3}{2}P - 2 \text{ si } P \equiv 0(\text{mod}2) \\ \frac{P}{2} + \frac{5}{2} \text{ si } P \equiv 0(\text{mod}2) \end{cases}$$

Les points particuliers de cette fonction sont :

(4.10) $$\begin{cases} P_{2,1} = 6 \\ P_{2,2} = 7 \\ A_2 = 5 \end{cases}$$

La figure suivante correspond à la matrice $\mathbb{T}_2(6,22,10)$

| | | | | | | | | | | |
|---|---|---|---|---|---|---|---|---|---|---|
| 6 | 7 | 6 | 7 | 6 | 7 | 6 | 7 | 6 | 7 | 6 |
| 7 | 6 | 7 | 6 | 7 | 6 | 7 | 6 | 7 | 6 | 7 |
| 8 | 10 | 13 | 9 | 7 | 6 | 7 | 6 | 7 | 6 | 7 |
| 9 | 7 | 6 | 7 | 6 | 7 | 6 | 7 | 6 | 7 | 6 |
| 10 | 13 | 9 | 7 | 6 | 7 | 6 | 7 | 6 | 7 | 6 |
| 11 | 8 | 10 | 13 | 9 | 7 | 6 | 7 | 6 | 7 | 6 |
| 12 | 16 | 22 | 31 | 18 | 25 | 15 | 10 | 13 | 9 | 7 |
| 13 | 9 | 7 | 6 | 9 | 7 | 6 | 7 | 6 | 7 | 6 |
| 14 | 19 | 12 | 16 | 22 | 31 | 18 | 25 | 15 | 10 | 13 |
| 15 | 10 | 13 | 9 | 7 | 6 | 7 | 6 | 7 | 6 | 7 |
| 16 | 22 | 31 | 18 | 25 | 15 | 10 | 13 | 9 | 7 | 6 |
| 17 | 11 | 8 | 10 | 13 | 9 | 7 | 6 | 7 | 6 | 7 |
| 18 | 25 | 15 | 10 | 13 | 9 | 7 | 6 | 7 | 6 | 7 |
| 19 | 12 | 16 | 22 | 31 | 18 | 25 | 15 | 10 | 13 | 9 |
| 20 | 28 | 40 | 58 | 85 | 45 | 25 | 15 | 10 | 13 | 9 |
| 21 | 13 | 9 | 7 | 6 | 7 | 6 | 7 | 6 | 7 | 6 |
| 22 | 31 | 18 | 25 | 15 | 10 | 13 | 9 | 7 | 6 | 7 |

Figure 17: La matrice $\mathbb{T}_2(6,22,11)$ avec le quatre couleurs

On peut vérifier aisément que :

(4.11) $$T^k(N) = F_4^k(N + 5) - 5$$

En se basant sur le théorème 2.1, la conjecture de Collatz peut être énoncée pour la fonction $F_2$ comme suit :

$\forall$ P $\in \mathbb{N}$ tel que P $\geq$ 6, il existe un entier naturel non nul k tel que :

$$F_2^k(P) = 6 \quad (4.12)$$

## 4.4 Comportements des suites générées par la fonction $F_{-1}$

La fonction $F_{-1}$ définie sur D(-1) est donnée par les expressions suivantes :

$$F_{-1}(P) = \begin{cases} \frac{3}{2}P + 1 & \text{si } P \equiv 0 \pmod{2} \\ \frac{1}{2}P - \frac{1}{2} & \text{si } P \equiv 1 \pmod{2} \end{cases} \quad (4.13)$$

Les points particuliers de cette fonction sont (n=-1):

$$\begin{cases} P_{-1,1} = 0 \\ P_{-1,2} = 1 \\ A_{-1} = -1 \end{cases} \quad (4.14)$$

On vérifie que :

$$F_{-1}(0) = 1 \text{ et } F_{-1}(1) = 0$$

Sur la figure suivante, on représente les deux matrices équivalentes $\mathbb{T}_{-1}(0,15,11)$ et $\mathbb{T}_{-\frac{1}{2}}(1,16,11)$

| | | | | | | | | | | |
|---|---|---|---|---|---|---|---|---|---|---|
| 0 | 1 | 0 | 1 | 0 | 1 | 0 | 1 | 0 | 1 | 0 |
| 1 | 0 | 1 | 0 | 1 | 0 | 1 | 0 | 1 | 0 | 1 |
| 2 | 4 | 7 | 3 | 1 | 0 | 1 | 0 | 1 | 0 | 1 |
| 3 | 1 | 0 | 1 | 0 | 1 | 0 | 1 | 0 | 1 | 0 |
| 4 | 7 | 3 | 1 | 0 | 1 | 0 | 1 | 0 | 1 | 0 |
| 5 | 2 | 4 | 7 | 3 | 1 | 0 | 1 | 0 | 1 | 0 |
| 6 | 10 | 16 | 25 | 12 | 19 | 9 | 4 | 7 | 3 | 1 |
| 7 | 3 | 1 | 0 | 1 | 0 | 1 | 0 | 1 | 0 | 1 |
| 8 | 13 | 6 | 10 | 16 | 25 | 12 | 19 | 9 | 4 | 7 |
| 9 | 4 | 7 | 3 | 1 | 0 | 1 | 0 | 1 | 0 | 1 |
| 10 | 16 | 25 | 12 | 19 | 9 | 4 | 7 | 3 | 1 | 0 |
| 11 | 5 | 2 | 4 | 7 | 3 | 1 | 0 | 1 | 0 | 1 |
| 12 | 19 | 9 | 4 | 7 | 3 | 1 | 0 | 1 | 0 | 1 |
| 13 | 6 | 10 | 16 | 25 | 12 | 19 | 9 | 4 | 7 | 3 |
| 14 | 22 | 34 | 52 | 79 | 39 | 19 | 9 | 4 | 7 | 3 |
| 15 | 7 | 3 | 1 | 0 | 1 | 0 | 1 | 0 | 1 | 0 |
| 1 | 2 | 1 | 2 | 1 | 2 | 1 | 2 | 1 | 2 | 1 |
| 2 | 1 | 2 | 1 | 2 | 1 | 2 | 1 | 2 | 1 | 2 |
| 3 | 5 | 8 | 4 | 2 | 1 | 2 | 1 | 2 | 1 | 2 |
| 4 | 2 | 1 | 2 | 1 | 2 | 1 | 2 | 1 | 2 | 1 |
| 5 | 8 | 4 | 2 | 1 | 2 | 1 | 2 | 1 | 2 | 1 |
| 6 | 3 | 5 | 8 | 4 | 2 | 1 | 2 | 1 | 2 | 1 |
| 7 | 11 | 17 | 26 | 13 | 20 | 10 | 5 | 8 | 4 | 2 |
| 8 | 4 | 2 | 1 | 2 | 1 | 2 | 1 | 2 | 1 | 2 |
| 9 | 14 | 7 | 11 | 17 | 26 | 13 | 20 | 10 | 5 | 8 |
| 10 | 5 | 8 | 4 | 2 | 1 | 2 | 1 | 2 | 1 | 2 |
| 11 | 17 | 26 | 13 | 20 | 10 | 5 | 8 | 4 | 2 | 1 |
| 12 | 6 | 3 | 5 | 8 | 4 | 2 | 1 | 2 | 1 | 2 |
| 13 | 20 | 10 | 5 | 8 | 4 | 2 | 1 | 2 | 1 | 2 |
| 14 | 7 | 11 | 17 | 26 | 13 | 20 | 10 | 5 | 8 | 4 |
| 15 | 23 | 35 | 53 | 80 | 40 | 20 | 10 | 5 | 8 | 4 |
| 16 | 8 | 4 | 2 | 1 | 2 | 1 | 2 | 1 | 2 | 1 |

Figure 18: Comportement des suites dans les deux matrices $\mathbb{T}_{-1}(0,15,11)$ et $\mathbb{T}_{-\frac{1}{2}}(1,16,11)$

Pour tous entiers naturels non nuls k et N, les deux fonctions $T^k$ et $F_{-1}^k$ vérifient l'équation suivante :

$$T^k(N) = F_{-1}^k(N-1) + 1 \quad (4.15)$$

En se basant sur le théorème 2.1, la conjecture de Collatz peut être énoncée pour la fonction $F_{-1}$ de la manière suivante :

$\forall$ P $\in \mathbb{N}$ (P $\geq$ 0), il existe un entier naturel non nul k tel que :

(4.16) $$\mathbf{F_{-1}^{k}(P) = 0}$$

## 4.5 Etude de la fonction $F_{-2}$

Cette fonction correspond à une valeur de paramètre n égale à -2, elle est définie sur D(-2) comme suit :

(4.17) $$F_{-2}(P) = \begin{cases} \frac{3}{2}P + 2 \text{ si } P \equiv 0 \ (\text{mod}2) \\ \frac{P}{2} - \frac{3}{2} \ \text{ si } P \equiv 1 \ (\text{mod}2) \end{cases}$$

Les points particuliers de cette fonction sont :

(4.18) $$\begin{cases} P_{-2,1} = -2 \\ P_{-2,2} = -1 \\ A_{-2} = -3 \end{cases}$$

On peut vérifier que :

$$F_{-2}(-2) = -1 \text{ et } F_{-2}(-1) = -2$$

La figure suivante correspond à la matrice $\mathbb{T}_{-2}(-2{,}13{,}11)$ noter que les cases sont colorées selon les règles présentées sur la figure 9.

| | | | | | | | | | | |
|---|---|---|---|---|---|---|---|---|---|---|
| -2 | -1 | -2 | -1 | -2 | -1 | -2 | -1 | -2 | -1 | -2 |
| -1 | -2 | -1 | -2 | -1 | -2 | -1 | -2 | -1 | -2 | -1 |
| 0 | 2 | 5 | 1 | -1 | -2 | -1 | -2 | -1 | -2 | -1 |
| 1 | -1 | -2 | -1 | -2 | -1 | -2 | -1 | -2 | -1 | -2 |
| 2 | 5 | 1 | -1 | -2 | -1 | -2 | -1 | -2 | -1 | -1 |
| 3 | 0 | 2 | 5 | 1 | -1 | -2 | -1 | -2 | -1 | -2 |
| 4 | 8 | 14 | 23 | 10 | 17 | 7 | 2 | 5 | 1 | -1 |
| 5 | 1 | -1 | -2 | -1 | -2 | -1 | -2 | -1 | -2 | -1 |
| 6 | 11 | 4 | 8 | 14 | 23 | 10 | 17 | 7 | 2 | 5 |
| 7 | 2 | 5 | 1 | -1 | -2 | -1 | -2 | -1 | -2 | -1 |
| 8 | 14 | 23 | 10 | 17 | 7 | 2 | 5 | 1 | -1 | -2 |
| 9 | 3 | 0 | 2 | 5 | 1 | -1 | -2 | -1 | -2 | -1 |
| 10 | 17 | 7 | 2 | 5 | 1 | -1 | -2 | -1 | -2 | -1 |
| 11 | 4 | 8 | 14 | 23 | 10 | 17 | 7 | 2 | 5 | 1 |
| 12 | 20 | 32 | 50 | 77 | 37 | 17 | 7 | 2 | 5 | 1 |
| 13 | 5 | 1 | -1 | -2 | -1 | -2 | -1 | -2 | -1 | -2 |

Figure 19: La Matrice $\mathbb{F}_{-2}(-2{,}13{,}11)$

En se basant sur le théorème 2.1, la conjecture de Collatz peut être énoncée pour la fonction $F_{-2}$ comme suit :

$\forall$ P $\in \mathbb{Z}$ tel que $P \geq -2$, il existe un entier naturel non nul k tel que :

$$F_{-2}^{k}(P) = -2 \tag{4.19}$$

## 4.6 Etude de la fonction $F_{-3}$

Cette fonction est définie sur D(-3) comme suit :

$$F_{-3}(P) = \begin{cases} \frac{3}{2}P + 3 \text{ si } P \equiv 0 \pmod 2 \\ \frac{P}{2} - \frac{5}{2} \text{ si } N \equiv 1 \pmod 2 \end{cases} \tag{4.20}$$

Les points particuliers de cette fonction sont :

$$\begin{cases} P_{-3,1} = -4 \\ P_{-3,2} = -3 \\ A_{-3} = -5 \end{cases} \tag{4.21}$$

On peut vérifier que :

$$F_{-3}(-4) = -3 \text{ et } F_{-3}(-3) = -4$$

La matrice $\mathbb{T}_{-3}(-4,11,10)$ est représentée sur la figure suivante les cases sont colorées selon les règles de l'équivalence chromatique données par la figure 9.

| | | | | | | | | | |
|---|---|---|---|---|---|---|---|---|---|
| -4 | -3 | -4 | -3 | -4 | -3 | -4 | -3 | -4 | -3 |
| -3 | -4 | -3 | -4 | -3 | -4 | -3 | -4 | -3 | -4 |
| -2 | 0 | 3 | -1 | -3 | -4 | -3 | -4 | -3 | -4 |
| -1 | -3 | -4 | -3 | -4 | -3 | -4 | -3 | -4 | -3 |
| 0 | 3 | -1 | -3 | -4 | -3 | -4 | -3 | -4 | -3 |
| 1 | -2 | 0 | 3 | -1 | -3 | -4 | -3 | -4 | -3 |
| 2 | 6 | 12 | 21 | 8 | 15 | 5 | 0 | 3 | -1 |
| 3 | -1 | -3 | -4 | -3 | -4 | -3 | -4 | -3 | -4 |
| 4 | 9 | 2 | 6 | 12 | 21 | 8 | 15 | 5 | 0 |
| 5 | 0 | 3 | -1 | -3 | -4 | -3 | -4 | -3 | -4 |
| 6 | 12 | 21 | 8 | 15 | 5 | 0 | 3 | -1 | -3 |
| 7 | 1 | -2 | 0 | 3 | -1 | -3 | -4 | -3 | -4 |
| 8 | 15 | 5 | 0 | 3 | -1 | -3 | -4 | -3 | -4 |
| 9 | 2 | 6 | 12 | 21 | 8 | 15 | 5 | 0 | 3 |
| 10 | 18 | 30 | 48 | 75 | 35 | 15 | 5 | 0 | 3 |
| 11 | 3 | -1 | -3 | -4 | -3 | -4 | -3 | -4 | -3 |

Figure 20: La Matrice $\mathbb{T}_{-3}(-4,11,10)$

En se basant sur le théorème 2.1, la conjecture de Collatz peut être énoncée pour la fonction $F_{-3}$ de la manière suivante :

$\forall$ P $\in \mathbb{Z}$ tel que $P \geq -4$, il existe un entier naturel non nul k tel que :

(4.22) $$F_{-3}^{k}(P) = -4$$

On fait comparer les deux matrices $\mathbb{T}_{-\frac{1}{2}}(1,16,10)$ et la matrice $\mathbb{T}_{-3}(-4,11,10)$ qui sont représentées sur la figure suivante :

| | | | | | | | | | |
|---|---|---|---|---|---|---|---|---|---|
| -4 | -3 | -4 | -3 | -4 | -3 | -4 | -3 | -4 | -3 |
| -3 | -4 | -3 | -4 | -3 | -4 | -3 | -4 | -3 | -4 |
| -2 | 0 | 3 | -1 | -3 | -4 | -3 | -4 | -3 | -4 |
| -1 | -3 | -4 | -3 | -4 | -3 | -4 | -3 | -4 | -3 |
| 0 | 3 | -1 | -3 | -4 | -3 | -4 | -3 | -4 | -3 |
| 1 | -2 | 0 | 3 | -1 | -3 | -4 | -3 | -4 | -3 |
| 2 | 6 | 12 | 21 | 8 | 15 | 5 | 0 | 3 | -1 |
| 3 | -1 | -3 | -4 | -3 | -4 | -3 | -4 | -3 | -4 |
| 4 | 9 | 2 | 6 | 12 | 21 | 8 | 15 | 5 | 0 |
| 5 | 0 | 3 | -1 | -3 | -4 | -3 | -4 | -3 | -4 |
| 6 | 12 | 21 | 8 | 15 | 5 | 0 | 3 | -1 | -3 |
| 7 | 1 | -2 | 0 | 3 | -1 | -3 | -4 | -3 | -4 |
| 8 | 15 | 5 | 0 | 3 | -1 | -3 | -4 | -3 | -4 |
| 9 | 2 | 6 | 12 | 21 | 8 | 15 | 5 | 0 | 3 |
| 10 | 18 | 30 | 48 | 75 | 35 | 15 | 5 | 0 | 3 |
| 11 | 3 | -1 | -3 | -4 | -3 | -4 | -3 | -4 | -3 |
| 1 | 2 | 1 | 2 | 1 | 2 | 1 | 2 | 1 | 2 |
| 2 | 1 | 2 | 1 | 2 | 1 | 2 | 1 | 2 | 1 |
| 3 | 5 | 8 | 4 | 2 | 1 | 2 | 1 | 2 | 1 |
| 4 | 2 | 1 | 2 | 1 | 2 | 1 | 2 | 1 | 2 |
| 5 | 8 | 4 | 2 | 1 | 2 | 1 | 2 | 1 | 2 |
| 6 | 3 | 5 | 8 | 4 | 2 | 1 | 2 | 1 | 2 |
| 7 | 11 | 17 | 26 | 13 | 20 | 10 | 5 | 8 | 4 |
| 8 | 4 | 2 | 1 | 2 | 1 | 2 | 1 | 2 | 1 |
| 9 | 14 | 7 | 11 | 17 | 26 | 13 | 20 | 10 | 5 |
| 10 | 5 | 8 | 4 | 2 | 1 | 2 | 1 | 2 | 1 |
| 11 | 17 | 26 | 13 | 20 | 10 | 5 | 8 | 4 | 2 |
| 12 | 6 | 3 | 5 | 8 | 4 | 2 | 1 | 2 | 1 |
| 13 | 20 | 10 | 5 | 8 | 4 | 2 | 1 | 2 | 1 |
| 14 | 7 | 11 | 17 | 26 | 13 | 20 | 10 | 5 | 8 |
| 15 | 23 | 35 | 53 | 80 | 40 | 20 | 10 | 5 | 8 |
| 16 | 8 | 4 | 2 | 1 | 2 | 1 | 2 | 1 | 2 |

Figure 21: La matrice $\mathbb{T}_{-\frac{1}{2}}(1,16,10)$ et la matrice $\mathbb{T}_{-3}(-4,11,10)$

Les distributions chromatiques pour toute les lignes de $\mathbb{T}_{-3}(-4,11,10)$ sont identiques à la distribution chromatiques relatives aux suites de Collatz dans la matrice $\mathbb{T}_{-\frac{1}{2}}(1,16,10)$.

Pour tous entiers naturel non nuls k et N, les deux fonctions $T^{k}$ et $F_{-3}^{k}$ vérifient l'équation suivante :

(4.23) $$T^{k}(N) = F_{-3}^{k}(N-5) + 5$$

## 4.7 Plongement des fonctions $F_n$ Extension des matrices pour les entiers relatifs $< 2n + 2$ :

Les deux théorèmes établies sont généraux, ils sont valables pour tous entiers relatifs du groupe des entiers relatifs $\mathbb{Z}$. Sur la figure suivante, on représente quelques suites dont les premiers termes sont des entiers relatifs inferieurs à 2n+2 dans ce cas on prend n=3 noté que les cases sont colorées selon les mêmes règles indiquées sur la figure 6 :

Cette fonction est définie comme suit :

(4.24)
$$F_3(P) = \begin{cases} \frac{3}{2}P - 3 \text{ si } P \equiv 0 \ (\text{mod}2) \\ \frac{P}{2} + \frac{7}{2} \quad \text{si } P \equiv 1 \ (\text{mod}2) \end{cases}$$

Les points particuliers de cette fonction sont :

(4.25)
$$\begin{cases} P_{3,1} = 8 \\ P_{3,2} = 9 \\ A_3 = 7 \end{cases}$$

Sur la figure suivante, on fait comparer les comportements des suites $S_{F,n}(P,k)$ pour des valeurs de $P < (2n+2)$ et les suites de Collatz $Sy(N,k)$ pour des valeurs de $N < 1$, noter que les cases sont colorées toujours selon les mêmes règles de la figure 9:

| | | | | | | | | | |
|---|---|---|---|---|---|---|---|---|---|
| 9 | 8 | 9 | 8 | 9 | 8 | 9 | 8 | 9 | 8 |
| 8 | 9 | 8 | 9 | 8 | 9 | 8 | 9 | 8 | 9 |
| 7 | 7 | 7 | 7 | 7 | 7 | 7 | 7 | 7 | 7 |
| 6 | 6 | 6 | 6 | 6 | 6 | 6 | 6 | 6 | 6 |
| 5 | 6 | 6 | 6 | 6 | 6 | 6 | 6 | 6 | 6 |
| 4 | 3 | 5 | 6 | 6 | 6 | 6 | 6 | 6 | 6 |
| 3 | 5 | 6 | 6 | 6 | 6 | 6 | 6 | 6 | 6 |
| 2 | 0 | -3 | 2 | 0 | -3 | 2 | 0 | -3 | 2 |
| 1 | 4 | 3 | 5 | 6 | 6 | 6 | 6 | 6 | 6 |
| 0 | -3 | 2 | 0 | -3 | 2 | 0 | -3 | 2 | 0 |
| -1 | 3 | 5 | 6 | 6 | 6 | 6 | 6 | 6 | 6 |
| 2 | 1 | 2 | 1 | 2 | 1 | 2 | 1 | 2 | 1 |
| 1 | 2 | 1 | 2 | 1 | 2 | 1 | 2 | 1 | 2 |
| 0 | 0 | 0 | 0 | 0 | 0 | 0 | 0 | 0 | 0 |
| -1 | -1 | -1 | -1 | -1 | -1 | -1 | -1 | -1 | -1 |
| -2 | -1 | -1 | -1 | -1 | -1 | -1 | -1 | -1 | -1 |
| -3 | -4 | -2 | -1 | -1 | -1 | -1 | -1 | -1 | -1 |
| -4 | -2 | -1 | -1 | -1 | -1 | -1 | -1 | -1 | -1 |
| -5 | -7 | -10 | -5 | -7 | -10 | -5 | -7 | -10 | -5 |
| -6 | -3 | -4 | -2 | -1 | -1 | -1 | -1 | -1 | -1 |
| -7 | -10 | -5 | -7 | -10 | -5 | -7 | -10 | -5 | -7 |
| -8 | -4 | -2 | -1 | -1 | -1 | -1 | -1 | -1 | -1 |

Figure 22: Comportements des suites générées par $F_3$ et des suites de Collatz dont les premiers termes sont inférieurs à (2n+2)

On remarque bien que dans le domaine $D^{-}(n)$ (sous ensemble de $\mathbb{Z}$ qui contient tous les entiers relatifs strictement inférieurs à 2n+2) les suites générées par la fonction $F_3$ se comportent d'une manière similaires aux comportements des suites de Collatz dans $\mathbb{Z}^{-}$.En réalité toutes les équations et les propriétés démontrées dans cet article sont valables pour touts les entiers relatifs sauf pour la conjecture de Collatz qui ne peut pas être énoncée pour une fonction $F_n$ que dans le domaine D(n).

## 5. Discussions des résultats

La conjecture de Collatz dans le cas d'une fonction $F_n$ énonce que pour tous entiers relatifs P supérieurs ou égaux à 2n+2, il existe un entier naturel non nul k tel que :

$$F_n^k(P) = 2n + 2$$

Selon la conjecture de Collatz on suppose que cette propriété est vraie pour entier relatif appartenant à D(n)

Si on fait tendre n vers$(-\infty)$, le domaine de vérification de cette conjecture tend vers $\mathbb{Z}$

(5.1) $$\lim_{n\to-\infty} D(n) = \mathbb{Z}$$

La fonction $F_n$ tend vers une fonction qu'on la note $F_{-\infty}$qui vérifie théoriquement la conjecture de Collatz pour tous entiers relatifs.

(5.2) $$\lim_{n\to+\infty} F_n = F_{-\infty}$$

Le domaine de la non vérification de la conjecture de Collatz tend vers $\mathbb{Z}$

(5.3) $$\lim_{n\to+\infty} D^{-}(n) = \mathbb{Z}$$

Le sous-ensemble $D^{-}(n)$ de $\mathbb{Z}$ est définie comme suit :

(5.4) $$D^{-}(n) = \left\{ P \in \mathbb{Z} \;\middle|\; P < 2n + 2 \right\}$$

Dans le cas contraire si on fait tendre n vers l'infini, on obtient théoriquement une fonction qui ne vérifie pas la conjecture de Collatz pour n'importe quel entier relatif.

(5.5) $$\lim_{n\to+\infty} P_{n,1} = \lim_{n\to+\infty} (2n + 2) = +\infty$$

Dans ce cas on obtient la fonction $F_{+\infty}$ tel que :

(5.6) $$\lim_{n\to-\infty} F_n = F_{+\infty}$$

Les résultats obtenus sont importants mais on doit pourtant développer ce travail, en se basant sur ces résultats afin de trouver une méthode pour prouver la conjecture de Collatz.

## Références